\documentclass[11pt,leqno]{article}

\usepackage{amsmath}
\usepackage{amssymb}
\usepackage{dsfont}
\usepackage{mathrsfs}
\usepackage{epsfig}
\usepackage{float}

\numberwithin{equation}{section}

\newtheorem{Theorem}{Theorem}

\newtheorem{Lemma}[Theorem]{Lemma}

\newtheorem{Definition}[Theorem]{Definition}
\newtheorem{remark}[Theorem]{Remark}

\numberwithin{Theorem}{section}

\bibstyle{plain}

\newcommand{\al}{\alpha}
\newcommand{\R}{{\mathbf R}}
\newcommand{\ds}{\displaystyle}
\newcommand{\e}{\varepsilon}

\newcommand{\lra}{\longrightarrow}
\newcommand{\ra}{\rightarrow}

\newcommand{\N}{{\mathbf N}}

\newcommand{\Q}{{\mathbf Q}}

\newcommand{\1}{{\mathds 1}}

\newcommand{\ph}{\varphi}

\newcommand{\Xo}{{\mathscr X}}

\newcommand{\les}{\leqslant}
\newcommand{\ges}{\geqslant}

\newcommand{\beq}{\begin{equation}}
\newcommand{\eeq}{\end{equation}}

\newcommand{\references}[1]{\theinstitutions
\begin{thebibliography}{#1}}

\newcommand{\Endrefs}{\end{thebibliography}}

\title{ Profile decomposition for sequences of  Borel measures}

\date{\it \small To the memory of Mircea Reghi\c s	}

\author{
Mihai MARI\c S\footnote{Institut de Math\'ematiques de
Toulouse UMR 5219, Universit\'e Paul Sabatier, 118 route de Narbonne, 31062 Toulouse
Cedex 9, France. {\sf e-mail}: mihai.maris@math.univ-toulouse.fr.}}

\begin{document}

\maketitle


\noindent
\begin{abstract}
We prove that, if dichotomy occurs when the concentration-compactness principle is used, 
the dichotomizing sequence can be choosen so that a nontrivial part of it concentrates. 
Iterating this argument leads to a profile decomposition for arbitrary sequences 
of bounded Borel measures. 
To illustrate our results we give an application to the structure of bouded sequences in the Sobolev space $ W^{1, p }( \R^N)$. 

\end{abstract}


\section{Introduction}

The concentration-compactness principle, introduced by P.-L. Lions in his celebrated papers \cite{lions}, \cite{lions-bis}
is an extremely effective tool in many areas in Mathematics and has  a huge number of applications. 
In its simplest form the principle asserts that given a sequence $(\mu_n)_{n \ges 1}$ of probability measures on $ \R^N$, 
there exists a subsequence $(\mu_{n_k})_{k \ges 1}$ that either spreads over the space so 
that the measure of each 
ball tends to zero uniformly with respect to the size of the ball
(this phenomenon is called "vanishing"; see Definition \ref{v} below for a precise statement), 
or $\mu_{n_k} $ can be "split" into two nontrivial parts whose supports are far away from each other (this situation is called "dichotomy"), 
or most  of the mass of $\mu_{n_k}$ remains on balls of fixed radius (the latter case is called "concentration"; for a precise statement see 
Definition \ref{c}). 

Many variants of the concentration-compactness principle exist in the literature 
(see, e.g., \cite{lions21}, \cite{lions22}; for "abstract" versions we refer to \cite{ST}, \cite{TF}). 
An important refinement is the so-called "profile decomposition" of sequences  of functions. 
Early variants of profile decomposition for sequences of approximate  solutions of some PDE 
can already be found in \cite{brezis-coron} and \cite{struwe}.
For arbitrary bounded sequences  in Sobolev spaces, 
profile decomposition was fully developed by P. G\'erard
(\cite{bahouri-gerard}, \cite{gerard},  \cite{GMO}). 
In the context of Strichartz estimates for the Schr\"odinger group $ e^{ i t \Delta }$, the first profile decomposition result was obtained by 
F. Merle and L. Vega \cite{merle-vega}.  
This technique  has led to important applications in the theory of nonlinear dispersive equations 
(see e.g. \cite{FV}, \cite{GG}, \cite{HK1}, \cite{HK2}, \cite{kenig-merle}, \cite{kenig-merle2}, \cite{keraani}, \cite{merle-vega}, \cite{raphael} and references therein). 

In most applications of the  concentration-compactness principle,  the difficult part is to understand the dichotomy case. 
We show in Lemma \ref{L6} that whenever dichotomy occurs, 
we can  choose a "dichotomizing subsequence"
$ \mu_{n_k} = \mu_{k,1} + \mu_{k, 2} + o(1)$ 
such that $ \mu_{k, 1}$ and $ \mu_{k, 2}$ have supports far away from each other and, in adition,  the sequence 
$( \mu_{k, 1})_{k \ges 1}$ "concentrates."
Iterating this argument we are able to prove a profile decomposition result for arbitrary sequences of bounded Borel measures. 
The spectrum of possible applications 
seems large, starting from the existence of minimizers in Calculus of Variations to 
possible extensions of the blowup theory for nonlinear dispersive equations initiated by C. Kenig and F. Merle (\cite{kenig-merle}, \cite{kenig-merle2}).  
For instance, Theorem \ref{main} below has already been crucial in proving the main results 
in \cite{MN}.

\medskip

\noindent
{\bf Definitions and notation.}
For the convenience of the reader we collect here the basic 
notions used throughout the paper.

\begin{Definition} \label{c}
Let $(\Xo, d)$ be a metric space and 
let $(\mu_n)_{n \ges 1}$ be a sequence of positive Borel measures on $ \Xo$.
Let $ (x_n)_{n \ges 1}$ be a sequence of points in $ \Xo$.
We say that $(\mu_n)_{n \ges 1}$ concentrates around $ (x_n)_{n \ges 1}$
if for any $ \e > 0$ there exists $ R_{\e } > 0 $ such that
$ \mu_n( \Xo \setminus B( x_n , R_{\e})) < \e $ for all $ n \ges 1$.
\end{Definition}

Sometimes we simply say that $(\mu_n)_{n \ges 1}$ concentrates if there is a sequence of points
around which $(\mu_n)_{n \ges 1}$ concentrates.
If $(\mu_n)_{n \ges 1}$  concentrates around $ (x_n)_{n \ges 1}$ it is obvious that any subsequence $(\mu_{ n_k})_{k \ges 1}$ concentrates
(around the corresponding subsequence $(x_{n_k})_{k \ges 1}$).

\begin{Definition} \label{v}
Let $ (\Xo,d)$ be a metric space and $(\mu _n)_{n \ges 1}$ a sequence of positive Borel measures on $\Xo$. 
We say that $(\mu _n)_{n \ges 1}$ 
is a vanishing sequence  if for any $r>0$, 
$$ \sup_{x \in \Xo} \mu_n (B(x,r)) \lra 0 \qquad \mbox{ as } n \lra \infty. $$
\end{Definition}

Given a positive Borel measure $ \mu$ on a metric space $\Xo$, the concentration function of $\mu $ is 
 $ q : [0, \infty) \lra [0, \infty)$ defined by $q(t) = \ds \sup_{x \in \Xo} \mu (B(x,t)) .$
It is clear that $q$ is nondecreasing and $ \ds \lim_{t \ra \infty } q(t) = \mu (\Xo).$
If $ (\mu_n )_{n \ges 1}$  is as above and $ q_n $ is the concentration function of $ \mu _n$, it is obvious that 
$ (\mu_n )_{n \ges 1}$  is a vanishing sequence if and only if 
$ \ds \limsup_{ n \ra \infty } q_n (t) = 0 $ for all $ t >0$, or equivalently $ \ds \lim_{ t \ra \infty} \left( \limsup_{ n \ra \infty } q_n (t) \right) = 0.$

Assume that $ \mu $ is a Borel measure on $ \Xo.$
If $ M \subset \Xo $ is a Borel set, we denote by $ \mu_{|M} $ the restriction of $ \mu $ to $M$  defined by $  \mu_{|M} (E) = \mu( M \cap E).$
Let $ (\rho_n)_{n \ges 1} \subset L^1( \Xo, \mu)$. 
Each $ \rho_n$ generates a Borel measure $ \mu_n$ defined by $ \mu_n(E) = \int_E \rho_n \, d \mu. $
The concentration function of $ \rho _n$ is the same as that of $ \mu _n$. 
We say that $ (\rho_n)_{n \ges 1}$ concentrates around $ (x_n)_{n \ges 1}$, respectively that it is a vanishing sequence, 
if $ (\mu_n)_{n \ges 1}$ has the same property.

\medskip

\noindent
{\bf Main results. }
Our main result is Theorem \ref{main} below. 
Its proof relies on the next lemma, which is of independent interest. 

\begin{Lemma}
\label{L6}
Let $ (\Xo, d)$ be a metric space and 
let $(\mu_n)_{n \ges 1} $ be a sequence of positive Borel measures on $ \Xo$ such that $\mu _n(\Xo) \les M$ for all $n$,
where $M$ is a positive constant.
Denote by  $ q_n(t) = \ds \sup_{ x \in \Xo} \mu_n (B(x, t)) $  the concentration function of $ \mu _n$.
Assume that 
\beq
\label{conv}
 q_n (t) \lra q(t) \quad \mbox{  as }  n \lra \infty  \mbox{  for a.e.  }  t >0  \qquad \mbox{ and  } \qquad \ds \lim_{t \ra \infty } q(t) = \al > 0. 
\eeq
Fix an increasing sequence $ t_n  \lra \infty $ such that $ q_n ( t_n)  \lra \al$    
 (the existence of such sequence is guaranteed by Lemma \ref{e} (ii)-(iii) below).

\medskip

There exist a subsequence $ (\mu_{n_k})_{k \ges 1}$, a sequence of points $ ( x_k)_{ k \ges 1} \subset \Xo$ and a sequence
$( R_k )_{ k \ges 1}$ such that $ R_k \les t_{n_k}$, $ R_k \lra \infty$, 
 $\mu_{n _k} ( B( x_k,  R_k )) \lra \al $ as $ k \lra \infty$,
and the sequence of measures $ ({\mu_{ n_k }}_{| B( x_k, R_k)})_{k \ges 1 } $ concentrates around $( x_k)_{k \ges 1}$. 

\end{Lemma}

\begin{remark} \rm 
If $ ({\mu_{ n_k }}_{| B( x_k, R_k)})_{k \ges 1 } $ concentrates around $( x_k)_{k \ges 1}$, it follows from Remark \ref{R2} (i) below that 
for any sequence $ t_k \lra \infty$, $ t_k < R_k$, the sequence  
 $ ({\mu_{ n_k }}_{| B( x_k, t_k)})_{k \ges 1 } $ also concentrates around $( x_k)_{k \ges 1}$ 
 and $ \mu_n ( B( x_k, R_k) \setminus B( x_k, t_k) )  \lra 0 $ as $ k \lra \infty$. 
In particular, for any $ \ph: [0, \infty ) \lra [0, \infty )$ such that $ \ph(t) < t $ and $\ph(t) \lra \infty$ as $ t \lra \infty$ we have 
$ \mu_{n_k} \left(  B( x_k, R_k) \setminus B( x_k, \ph( R_k) ) \right)  \lra 0 $ as $ k \lra \infty$.
In applications it is very important to dispose of a "large" region $ B( x_k, R_k) \setminus B( x_k, \ph( R_k) )  $ with "small" measure
in order to perform a convenient cut-off (see also Remark \ref{R1.7} below).

\end{remark}

\begin{Theorem} \label{main}
Let $(\Xo,d) $ be a metric space and $(\mu_n)_{n \ges 1}$ a sequence of positive Borel measures on $\Xo$ such that 
$$
M:= \limsup_{n \ra \infty} \mu_n (\Xo) < \infty. 
$$
Let  $ \ph : [0, \infty) \to [0, \infty) $ be an increasing function such that $ \ph (s) \! \les \! \frac s2$ for all $s$ and 
$ \ds \lim_{s \ra \infty} \ph(s) \! = \infty $.

Then either $ ( \mu _n)_{n \ges 1}$ is a vanishing sequence, or 
there exists an increasing mapping $ j : \N^* \lra \N^* $ such that the subsequence $ (\mu_{j(n)} )_{n \ges 1} $ satisfies 
one of the following properties:



\medskip

i) There are $k \in \N^*$, positive numbers $ m_1, \dots m_k$, sequences of points $(x_n^i)_{n \ges 1} \subset \Xo$ 
and increasing sequences of positive numbers $(r_n^i)_{n \ges 1}$ such that $ r_n ^i \lra \infty $ as $ n \lra \infty$, 
$ i \in \{ 1, \dots , k \}$ satisfying the following properties:

\medskip

$\qquad $ a) For each $n$ the balls $B(x_n ^i, r_n^i)$, $ i \in \{ 1, \dots , k \}$ are disjoint.

\medskip

$\qquad $ b) For each $ i \in \{ 1, \dots , k \}$ we have 
$$
\mu_{ j(n)} \left( B( x_n ^i, \ph ( r_n ^i)) \right) \lra m_i \qquad \mbox{ as } n \lra \infty ,  \qquad 
\mu_{ j(n)} \left( B( x_n ^i,  r_n ^i )  \setminus B( x_n ^i, \ph ( r_n ^i)) \right) \les \frac{1}{2^{n+i}}.
$$
and the sequence of measures $\left( { \mu_{ j(n)} }_{ | B( x_n ^i,  r_n ^i ) } \right) _{n \ges 1} $ 
concentrates around $(x_n^i)_{n \ges 1}$.

\medskip

$\qquad $ c) The sequence of measures 
$\left( { \mu_{ j(n)} }_{ \big| \Xo \setminus \cup_{i =1}^k B( x_n ^i,  r_n ^i ) } \right) _{n \ges 1} $  is a vanishing sequence.

\medskip

ii) There are  positive numbers $ m_1, \dots m_k, \dots $ such that $ m_{k+1} \les 2 m_k$,  
sequences of points $(x_n^i)_{n \ges i} \subset \Xo $ 
and increasing sequences of positive numbers $(r_n^i)_{n \ges i}$ such that $ r_n ^k \lra \infty $ as $ n \lra \infty$ for each fixed $k$ 
and the following properties hold:

\medskip

$\qquad $ a) For each $n$ the balls $B(x_n ^1, r_n^1), \dots , B( x_n ^n, r_n ^n) $   are disjoint.

\medskip

$\qquad $ b) The same as  (b) in (i) above. 

\medskip

$\qquad $ c) 
Denote by   $ \tilde{q} _n ^{\ell } $ is the concentration function of $ { \mu_{ j(n)} }_{ \big| \Xo \setminus \cup_{i =1}^{\ell} B( x_n ^i,  r_n ^i ) }$
for  $ \ell \ges n.$  Then 
$$
\ds \lim_{\ell \ra \infty} \left( \lim_{t \ra \infty} \left( \limsup_{n \ra \infty} \tilde{q} _n ^{\ell } (t) \right) \right) = 0. 
$$

$\qquad $ d) The sequence of measures 
$\left( { \mu_{ j(n)} }_{ \big| \Xo \setminus \cup_{i =1}^n B( x_n ^i,  r_n ^i ) } \right) _{n \ges 1} $   is a vanishing sequence.

\end{Theorem}

\begin{remark} 
\rm
\label{R1.7}
Let $(\ph _k)_{k \ges 1}$ be a sequence of increasing functions, where each $ \ph _k$ is as in Theorem \ref{main} and
$ \ph_1 \ges \ph _2 \ges \dots \ges \ph_k \ges \ph_{k+1} \ges \dots $. 
Theorem \ref{main} above still holds true (with obvious modifications in the proof) if  (ii) (b) is replaced by  
$$
\mu_{ j(n)} \left( B( x_n ^k, \ph _k ( r_n ^k)) \right) \lra m_k \quad \mbox{ as } n \lra \infty ,  \qquad 
\mu_{ j(n)} \left( B( x_n ^k,  r_n ^k )  \setminus B( x_n ^k, \ph _k ( r_n ^k)) \right) \les \frac{1}{2^{n+k}}.
$$
The region $B( x_n ^k,  r_n ^k )  \setminus B( x_n ^k, \ph _k ( r_n ^k)) $ has small measure and, in applications, 
is often used to perform an appropriate "cut-off."
The fact that the $\ph _k$'s may depend on $k$ and may tend to infinity as slowly as we want 
 gives even more freedom to perform the cut-off. 
 We refer to \cite{MN} for an application ($\ph(t) = t^{\frac 13}$ has been used there). 
\end{remark}

\medskip

As an application 
we prove a "profile decomposition" result for bounded sequences in the Sobolev space $W^{1, p } ( \R^N)$. 
Our proof is simple, direct and relies only on Theorem \ref{main}.
Similar results already exist in the literature. 
To our knowledge, one of the first results of this kind was Theorem 1.1 p. 214 in \cite{gerard} which deals with bounded sequences in $\dot{H}^s( \R^N)$, $ 0< s < \frac N2.$
The proofs in \cite{gerard} use in an essential way the Fourier transform and improved Sobolev inequalities. 
An extension to $\dot{H}^{s, p}( \R^N)$ spaces is given by Theorem 1 p. 386 in \cite{jaffard}. 
The proof in \cite{jaffard} relies on wavelet decomposition and improved Sobolev inequalities. 
We have to mention the recent work \cite{PP}, which recovers the profile decomposition in \cite{gerard} in $ \dot{H}^s $ 
(among other interesting results) and uses again improved Sobolev estimates.  
The folowing result is a generalization to $ W^{1, p}$ spaces of 
Proposition 3.1 p. 2822 in \cite{HK1} and of Proposition 2.1 p. 1039 in \cite{HK2}, which deal with bounded sequences in $H^1( \R^N)$. 
Unlike in the previous works we do not use the Fourier transform and improved Sobolev inequalities, nor the Hilbert space structure. 
We postpone to a subsequent work a similar study of bounded sequences in the spaces $W^{s, p} (\R^N)$ 
as well as in the spaces
$\dot{W} ^{1, p }( \R^N )$ (the latter case is a little bit more delicate because of scaling invariance, but Theorem \ref{main} can still be used to get a profile decomposition result 
in the spirit of \cite{gerard}, \cite{bahouri-gerard}).

\begin{Theorem}
\label{Sobolev}
Let $ 1 < p < \infty$ and let $(u_n)_{n \ges 1}$ be a bounded sequence in $ W^{1, p }(\R^N)$. 
Let $ p^* = \frac{ Np}{N-p}$ if $ p < N$ and  $ p^ * = \infty $ if $ p \ges N$. 
There exists a subsequence of $(u_n)_{n \ges 1}$ (still denoted the same), 
a family $( V^i)_{i \in \N^*}$ of functions in $ W^{1, p }(\R^N)$ and a family of sequences $( y_n^i)_{n \ges 1} \subset \R^N$, $i \in \N^*$, such that: 

\medskip

i) For any $ i \neq j$ we have $ | y_n ^i - y_n ^j| \lra \infty $ as $ n \lra \infty$. 

\medskip

ii) For any  $ k \in \N^*$ there holds 
$$
 u_n = \ds \sum_{i = 1}^k V^i ( \cdot - y_n ^i) + w_n ^k , \qquad
\mbox{ where } \ds \lim_{ k \ra \infty } \left( \limsup_{n \ra \infty} \| w_n ^k \|_{ L^q( \R^N)} \right) = 0 
\mbox{ for all } q \in (p, p^*). 
$$


iii) For all  $ k \ges 1$ we have $ \| u_n \|_{L^p ( \R^N)}^p = \ds \sum_{i =1}^k \| V^i \|_{L^p (\R^N)} ^p + \| w_n ^k  \|_{L^p ( \R^N)}^p + o(1) $ as $ n \lra \infty$. 

iv)
If $ p = 2$ we have 
$ \| \nabla u_n \|_{L^2 ( \R^N)}^2 \! = \! \ds \sum_{i =1}^k \| \nabla  V^i \|_{L^2 (\R^N)} ^2 \! +  \| \nabla w_n ^k  \|_{L^2 ( \R^N)}^2 + o(1) $ as $ n \! \lra \! \infty$. 
\end{Theorem}

As pointed out in \cite{jaffard} p. 387,  the equality in Theorem \ref{Sobolev} (iv) may be false  if $ p \neq 2$. 

\medskip

The rest of this paper is organized as follows.
In the next section we give some elementary results that will be used in proofs. 
We prove Lemma \ref{L6} in section \ref{PrL6} and Theorem \ref{main} in section \ref{Prmain}.
The proof of Theorem \ref{Sobolev} is given in section \ref{PrSobolev}.

\section{Some elementary facts}

In this section we collect some simple observations that will be very useful in the sequel. 

\begin{remark} \label{R2} \rm Let $(\mu_n)_{n \ges 1}$ be a sequence of positive Borel measures on $ \Xo$.

\medskip

(i)  Assume that  $(\mu_n)_{n \ges 1}$ concentrates around $ (x_n)_{n \ges 1}$.  If $(M_n)_{n \ges 1} $ are any Borel measurable sets,
the sequence $({\mu_n}_{\big| M_n} )_{n \ges 1}$ also concentrates around $ (x_n)_{n \ges 1}$.
If $ (t_n)_{n \ges 1} $ is any sequence of positive numbers such that $ t_n \lra \infty$, then
$\mu_n ( \Xo \setminus B( x_n, t_n)) \lra 0 $ as $ n \lra \infty$.

\medskip

(ii)  If $(\mu_n)_{n \ges 1}$ concentrates around $ (x_n)_{n \ges 1}$ and $(y_n)_{n \ges 1}$ is any sequence of points such that 
$(d( x_n ,  y_n))_{n \ges 1}$  is bounded, then  $(\mu_n)_{n \ges 1}$ concentrates around $ (y_n)_{n \ges 1}$.

\medskip

(iii) Denote $ \al _n = \mu_n ( \Xo)$. Assume that each $ \al _n $ is finite and $ \al _n \ges \al _0$, where $ \al _0 $ is a positive constant.
Let $ q_n ( r ) = \ds \sup _{x \in \Xo} \mu_n( B(x, r)) $ be the concentration function of $ \mu_n$.
Then $(\mu_n)_{n \ges 1}$ concentrates if and only if
\beq
\label{1}
\forall \e > 0, \; \exists r_{\e } > 0, \; \forall r \ges r_{\e }, \; \forall n \ges 1, \; \quad 0 \les \al _n - q_n (r) < \e.
\eeq

{\it Proof. } If $(\mu_n)_{n \ges 1}$ concentrates around $ (x_n)_{n \ges 1}$ and $ R_{\e }$ is as in Definition \ref{c},
for any $ r \ges R_{\e}$ we have
$
\al _n \ges q_n(r) \ges q_n (R_{\e} ) \ges \mu_n ( B( x_n, R_{\e})) > \al _n - \e
$
and (\ref{1}) holds.

Conversely, assume that (\ref{1}) is satisfied and fix $ \e \in (0, \frac{ \al _0}{2}).$
Let $ r_{\e}$ be given by (\ref{1}). For all $ n$ we have $q_n ( r_{\e }) > \al _n - \e$,
hence we may choose  $ x_n \in \Xo$ such that $ \mu_n ( B( x_n, r_{\e}) ) > \al _n - \e$.
Let $ 0< \e'< \e$.
Choose $ r_{\e'}$ such that (\ref{1}) holds.
In particular we have $ q_n ( r_{\e'}) > \al _n - \e'$ for all $n$, hence there exist $ y_{n, \e'} \in \Xo$ such that
$ \mu_n (B(y_{n, \e'}, r_{\e'}) ) > \al _n - \e'.$
The balls $B( x_n, \e) $ and $B(y_{n, \e'}, r_{\e'})  $ must intersect (for otherwise we would have
$ \al_n = \mu_n ( \Xo) \ges \mu_n( B( x_n, \e)) + \mu_n( B(y_{n, \e'}, r_{\e'}) ) 
> \al _n - \e + \al_n - \e '>  \al _n  $,
a contradiction). Hence $ B(y_{n, \e'}, r_{\e'}) \subset B( x_n, r_{\e} + 2 r_{\e'}).$
Taking $ R_{\e'} = r_{\e} + 2 r_{\e'}$ for $ 0 < \e'< \e$ we have
$ \mu _n ( \Xo \setminus B(x_n, R_{\e '})) \les  \mu _n ( \Xo \setminus  B(y_{n, \e'}, r_{\e'}) ) < \e '$, 
therefore  $(\mu_n)_{n \ges 1}$ concentrates around $ (x_n)_{n \ges 1}$.

\medskip

(iv) We keep the previous notation and we assume that  $ (\al _n)_{n \ges 1} $ is bounded,
$ 0 < \al _0 \les \al _n$ for all $ n \ges 1$
and there exists an increasing sequence
$ t_k \lra \infty$ such that  for any $k$ the sequence $(q_n (t_k))_{n \ges 1}$ converges.
Denote $ q( t_k) = \ds \lim_{n \ra \infty} q_n( t_k) $ and $ Q = \ds \lim _{k \ra \infty } q( t_k) $
($Q$ exists because $(q(t_k))_{k \ges 1}$ is nondecreasing).
Then we have  $ Q \les \ds \liminf _{ n \ra \infty} \al _n$. Moreover, 
 $ (\mu_n)_{n \ges 1}$ concentrates if and only if $ \ds \lim_{n \ra \infty} \al _n = Q$.

{\it Proof. }
 We have $ q_n ( t_k) \les \al _n $ for all $ n $ and $ k$, hence $ q( t_k) \les  \ds \liminf _{ n \ra \infty} \al _n$ for all $k$
and the first claim follows.

Assume that $ (\mu_n)_{n \ges 1}$ concentrates.
Fix $ \e >0$.
Let $ r_{\e}  $ be as in (\ref{1}).
There is $ k_{\e } $ such that $ t_{k_{\e}} > r_{\e } $ and $ | q(t_{k_{\e}}) - Q | < \e$.
By (\ref{1}) we have $ \al _n - \e < q_n( t_k) \les \al _n$ for all $ k \ges k_{\e}$ and $ n \ges 1$.
Choose $ n_{\e}$ such that $ | q_n ( t_{k_{\e}}) - q(t_{k_{\e}}) | < \e$ for $ n \ges n_{\e}$.
Then
$$
| \al _n - Q| \les | \al _n - q_n (t_{k_{\e}}) | + |q_n (t_{k_{\e}})  - q(t_{k_{\e}}) |  + |  q(t_{k_{\e}})  - Q| < 3 \e
\qquad \mbox{ for all } n \ges n_{\e}.
$$
Since $ \e$ is arbitrary we infer that $ \al _n \lra Q$ as $ n \lra \infty$.

Conversely, assume that  $ \al _n \lra Q$.
Fix $ \e > 0$. Choose $ k_{\e}$ such that $|  q( t_{k_{\e}}) - Q| < \frac{\e}{3}$, then choose $ n_{\e}$ such that
$ | \al _n - Q| < \frac{ \e}{3}$ and $ | q_n (t_{k_{\e}}) -  q( t_{k_{\e}}) | < \frac{\e}{3} $ for all $ n \ges n_{\e}$.
Since $ q_n $ is nondecreasing, for all $ t \ges t_{k_{\e}} $ and $ n \ges n_{\e}$ we have
$$
0 \les \al _n - q_n(t) \les \al_n - q_n (t_{k_{\e}})
 \les | \al _n - Q| + |Q -  q( t_{k_{\e}}) | + |  q( t_{k_{\e}}) -  q_n( t_{k_{\e}}) | < \e.
$$
Since $ q_n (t) \lra \al _n$ as $ t \lra \infty$ for any fixed $n$, 
we may chose $ r_{\e} \ges t_{k_{\e}} $ such that
$ 0 \les   \al _n - q_n(t) < \e$ for all $n \in \{ 1, \dots, n_{\e} -1 \} $ and $ t \ges r_{\e}$.
Hence (\ref{1}) holds, that is $ (\mu_n)_{n \ges 1}$ concentrates.

\end{remark}

We will repeatedly use the following  lemmas.

\begin{Lemma} \label{e}
Let $ (f_n)_{n \ges 1}$ be a sequence of nondecreasing functions, $ f_n : [0, \infty ) \lra [0, \infty)$.
Assume that there is an increasing sequence $ (s_k)_{k \ges 1} \subset [0, \infty) $ such that
$ s_k \lra \infty$ 
and for each $k$ the sequence $(f_n(s_k))_{n \ges 1}$ has a limit that we denote  $ f( s_k)$.
Let
 $ \ell = \ds \lim_{k \ra \infty} f( s_k) $ (the limit exists because $(f( s_k))_{k \ges 1}$ is nondecreasing). Then:

\medskip

(i) For any sequence $ t_n \lra \infty $ we have $ \ds \liminf _{n \ra \infty} f_n ( t_n) \ges \ell. $

\medskip

(ii) There exists a nondecreasing sequence $ (t_n^*)_{n \ges 1}$, $ t_n ^* \lra \infty$ such that $ f_n ( t_n ^* ) \lra \ell.$

\medskip

(iii) If  $ (t_n^*)_{n \ges 1}$ is as in (ii), for any sequence $ (t_n)_{n \ges 1}$ such that $ 0 \les t_n \les t_n ^* $ and $ t_n \lra \infty$
we have $ f_n ( t_n) \lra \ell. $

\end{Lemma}

{\it Proof. } (i) Fix $ \ell '< \ell. $ Choose $ k$ such that $ f( s_k) > \ell'$.
There is $ n_0 $ such that $ f_n ( s_k) > \ell '$ for all $ n \ges n_0$.
There is $ n_1 \ges n_0$ such that $ t_n \ges s_k$ for all $ n \ges n_1$.
Then for any $ n \ges n_1$ we have $ f_n ( t_n) \ges f_n ( s_k) > \ell '$ and we infer that 
$ \ds \liminf_{n \ra \infty} f_n ( t_n) \ges \ell'.$ Since $ \ell'< \ell $ is arbitrary, 
(i) follows.

\medskip

(ii) 
If $ \ell = \infty$, it follows from (i) that any sequence $ t_n \lra \infty$ satisfies (ii).
Assume next that $ \ell $ is finite.
Then $ 0 \les f(s_k) < \infty $ for all $k$.
We choose inductively $n_1, \dots , n_k , \dots \in \N$ such that
$ n_k < n_{k +1} $ and
$$
| f_n (s_k) - f( s_k) | < \frac 1k \qquad \mbox{ for all } n \ges n_k.
$$
If $ n_k \les n < n_{k+1}$ put $ t_n ^* = s_k$.
Then for $ n_k \les n < n_{k+1}$ we have
$$
| f_n ( t_n^*) - \ell | = | f_n ( s_k) - \ell | \les
| f_n ( s_k) - f( s_k) | + | f( s_k) - \ell | < \frac 1k + | f( s_k) - \ell | \lra 0
$$
as $ k \lra \infty$, hence $ t_n ^* $ satisfies (ii).

\medskip

Assertion (iii) is a simple consequence of (i) and (ii).

\hfill
$\Box$

\begin{Lemma} \label{e2}
Let $ (q_n)_{n \ges 1}$, $ q_n : [0, \infty) \lra \R$ be a sequence of nondecreasing functions such that
$(q_n(t))_{n \ges 1}$ is bounded
for any fixed $ t \ges 0$.
Denote $ \tilde{q}(t) = \ds \limsup_{ n \ra \infty} q_n(t)$ and $ \al = \ds \lim_{ t \ra \infty} \tilde{q}(t) $
(the last limit exists because $ \tilde{q}$ is nondecreasing).

For any $ \beta < \al $ there is a subsequence $ (q_{n_k})_{ k \ges 1}$ and a nondecreasing function $ q : [0, \infty) \lra \R$
such that 
$ q_{n_k}(t) \lra q(t) $ as $ k \lra \infty$ for any $t\ges 0$ and $ \ds \lim_{ t \ra \infty} {q}(t) > \beta .$
\end{Lemma}

{\it Proof. } Fix $ \beta < \al$.
Then fix $ t_0 > 0 $ such that $ \tilde{q} ( t_0 ) > \beta. $
There is an increasing mapping $ j_0 : \N ^* \lra \N^* $ such that $ q_{j_0 (n) } (t_0) \lra \tilde{q} (t_0) > \beta$.
We may write $ \Q_+ \setminus \{ t_0 \} = \{ t_1, t_2, \dots, t_n , \dots \}$, where  $ t_i \neq t_j $  if $ i \neq j$.

Inductively we construct increasing mappings $ j_1, j_2, \dots, j_n , \dots \; : \N^* \lra \N^*$ such that
$ q_{ j_0( j_1(n))} ( t_1)$ converges to a limit denoted $ q( t_1) $,
$ q_{ j_0( j_1(j_2(n)))} ( t_2)$ converges to a limit denoted $ q( t_2) $ and so on,
$ q_{ j_0( j_1( j_2( \dots j_k(n) \dots )))} ( t_k)$ converges as $ n \lra \infty $ to a limit denoted $ q( t_k) $
for any $ k \ges 1 $.
Let $ j (n) = j_0( j_1( \dots ( j_n(n)) \dots ))$ and $ q( t_0) = \tilde{ q} ( t_0)$.
It is clear that $ j: \N^* \lra \N^*$ is increasing,
$ q_{j(n)}(t) \lra q(t) $  as $ n \lra \infty $ for all $ t \in \Q_+ \cup \{ t _0\}$
and $ q $ is a nondecreasing mapping defined on $ \Q_+ \cup \{ t _0\}$.

We extend $ q$ to a nondecreasing mapping on $ [0, \infty)$ by  defining
$ q(t) = \ds \sup \{ q(s) \; | \; s \in \Q_+ \cup \{ t_0 \} , \; s \les t \} $.
We claim that $ q_{j(n)}(t) \lra q(t)$ as $ n \lra \infty$ at any $ t > 0 $ where $q$ is continuous.
Fix $ t > 0 $ such that $ q $ is continuous at $ t$. Fix $ \e > 0$.
There are $ r, s \in \Q_+ $ such that $ r < t < s$ and
$ q( t ) - \frac { \e}{2} < q(r) \les q(t) \les q(s) < q( t) + \frac{ \e}{2}.$
There is $ n_{\e} \in \N $ such that $ q_{j(n)} ( r) > q(r) - \frac{ \e}{2} $ and $ q_{j(n)} ( s) < q(s) + \frac{ \e}{2} $
for all $ n \ges n_{\e}$.
Hence
$$
q(t) - \e <  q(r) - \frac{ \e}{2} <  q_{j(n)} ( r) \les  q_{j(n)} ( t) \les  q_{j(n)} ( s) < q(s) + \frac{ \e}{2} < q(t) + \e
$$
for all $ n \ges n_{\e}$.
Since $ \e $ is arbitrary this implies $  q_{j(n)} ( t) \lra q(t)$.
Moreover, $q(t_0) = \tilde{q} ( t_0) > \beta$, hence $ q(t) \ges q( t_0) > \beta $ for all $ t \ges t_0$.

There are at most countably many points  where $q$ is discontinuous, say $ d_1, d_2, \dots, d_k, \dots$
Using again the diagonal extraction procedure we construct an increasing mapping $ \kappa : \N^* \lra \N^*$ 
such that $ (q _{j( \kappa (n))} (d_{\ell})) _{n \ges 1}$ converges for any $ \ell \ges 1$. 
Then we redefine $ q( d _{\ell} ) =\ds  \lim_{n \ra \infty} q _{j( \kappa (n))} (d_{\ell}).$
It is obvious that $q$ and the subsequence $ (q _{j( \kappa (n))}  )_{n \ges 1}$ satisfy the conclusion of Lemma \ref{e2}
\hfill
$\Box$

\begin{remark} \rm
Lemma \ref{e2} is a refinement of Helly's Lemma. 
It is optimal in the following sense:
we may construct a sequence of nondecreasing functions $ q_n : [0, \infty) \lra [0, 1]$ such that
$\ds \lim _{t \ra \infty} \left( \limsup_{ n \ra \infty} q_n(t) \right) = 1 $
and for any subsequence $(q_{n_k})_{k \ges 1}$ such that there exists  $ t_0 > 0$ satisfying
$ \ds \lim_{k \ra \infty } q_{n _k} ( t_0 ) > 0,$
there holds $ \ds \lim _{t \ra \infty} \left( \limsup_{ k \ra \infty} q_{n_k}(t) \right) < 1 .$

\end{remark}

\section{Proof of Lemma \ref{L6}} 
\label{PrL6}

Fix an increasing function $ \ph \! : \! [0, \infty ) \! \lra \! [0, \infty)$ such that 
$ \ph (t) \les \! \frac t2 $ for all $t$ and $ \ds \lim_{ t \ra \infty } \ph(t) \!  = \! \infty$.
We prove first a relaxed version of Lemma \ref{L6}  which  asserts only that 
$ ({\mu_{ n_k }}_{| B( x_k, R_k)})_{k \ges 1 } $ concentrates (not necessarily around $( x_k)_{k \ges 1}$)
and 
\beq
\label{annulus}
\mu_{n _k} \left( B( x_k, R_k)  \setminus B( x_k, \ph ( R_k)) \right) \lra 0 \qquad \mbox{  as }  k \lra \infty.
\eeq
It follows from Lemma \ref{e} (iii) that $ q_n ( \ph  ( t_n)) \lra \al$.
For each $n$ there is $ x_n ^0 \in \Xo$ such that $ \mu_n ( B( x_n^0, \ph ( t_n))) > q_n ( \ph( t_n) ) - \frac 1n.$
It is clear that
$ q_n (t_n ) \ges \mu _n ( B( x_n^0, t_n)) \ges \mu _n ( B( x_n^0, \ph(t_n))) $, therefore
\beq
\label{2}
\mu _n \left( B( x_n^0, \ph(t_n))\right)  \lra \al \quad \mbox{ and } \quad
\mu _n \left( B(x_n^0, t_n ) \setminus B( x_n^0, \ph(t_n)) \right) \lra 0
\quad \mbox{ as } n \lra \infty.
\eeq

Denote by $ q_n ^1 (t)$ the concentration function of $ {\mu _n}_ {| B( x_n ^0, t_n )}.$
For each $ t \ges 0$ let $ \tilde{q}^1(t) = \ds \limsup_{ n \ra \infty	} q_n^1(t). $
Using Lemma \ref{e2} we see that there are a nondecreasing function $ q^1 : [0, \infty) \lra [0, M]$
and an increasing mapping  $ j_1 : \N^* \lra \N^*$ such that
\beq
\label{3}
q_{ j_1 (n)} ^1 (t) \lra q^1 (t) \qquad \mbox{ for all  } t > 0 \mbox{ and }
\lim_ { t \ra \infty} q^1(t) > \lim_ { t \ra \infty} \tilde{q} ^1(t) - \frac{ \al }{2^2}.
\eeq

Let $ \al _1 = \ds \lim_ { t \ra \infty} q ^1(t).$
Since $ q_n ^1 \les q_n$, it is clear that $ q^1 \les q $ and $ 0 \les \al _1 \les \al.$
If $ \al _1 = \al $ it follows from (\ref{2}) and Remark \ref{R2} (iv) that
$\left( {\mu _{ j_1(n)}}_ {| B( x_{j_1(n)} ^0, t_{j_1(n)} )}\right)_{n \ges 1} $
concentrates,
hence the conclusion of Lemma \ref{L6} (relaxed) is satisfied by 
$\left( \mu_{j_1(n)} \right)_{n \ges 1}$,
$\left( x_{j_1(n)}^0 \right)_{n \ges 1}$ and   $\left( t_{j_1(n)} \right)_{n \ges 1}$.

Assume that $ \al _1 < \al.$
By Lemma \ref{e} (ii)-(iii)  there is an increasing  sequence $ t_n ^1 \lra \infty $
such that $ t_n ^1 < \frac 14 t_{j_1(n)} $ and $ q_{j_1(n)}^1 ( t_n ^1) \lra \al _1$.
Then using the fact that $ q_{j_1(n)} ( t_{j_1(n)}) \lra \al $, (\ref{conv}) and
Lemma \ref{e} (iii) we get 
$$
q_{j_1(n)} ( t_n^1) \lra \al \qquad \mbox{ and } \qquad q_{j_1(n)} ( \ph( t_n^1))  \lra \al \qquad \mbox{ as } n \lra \infty.
$$
For each $n$ there is $ x_n ^1 \in \Xo$ such that
\beq
\label{4}
\mu_{ j_1 (n)} \left( B ( x_n ^1, \ph ( t_n ^1)) \right) > q_{j_1(n)} \left( \ph( t_n^1) \right) - \frac{1}{ 2 j_1 (n)} .
\eeq
It is clear that
$ \mu_{ j_1 (n)} \left( B ( x_n ^1, \ph ( t_n ^1)) \right) \les \mu_{ j_1 (n)} \left( B ( x_n ^1,  t_n ^1) \right) \les q_{j_1(n)} ( t_n ^1), $
hence
\beq
\label{5}
\mu_{ j_1 (n)} \left( B ( x_n ^1, \ph ( t_n ^1)) \right) \lra \al \qquad \mbox{ and } \qquad
\mu_{ j_1 (n)} \left( B ( x_n ^1,  t_n ^1 ) \right) \lra \al \qquad \mbox{ as } n \lra \infty .
\eeq

We claim that there is $ n_1 \in \N^*$ such that for all $ n \ges n_1$,
\beq
\label{6}
B ( x_{j_1(n)}^0, \ph( t_{j_1(n)}))  \cap B( x_n ^1, t_n ^1) = \emptyset.
\eeq
To see this we argue by contradiction and assume that the balls  intersect.
Recalling that  $ \ph( t_{j_1(n)}) \les \frac 12 t_{j_1(n)} $ and $ t_n ^1 < \frac 14 t_{j_1(n)} $, we have 
$$
B( x_n ^1, t_n ^1) \subset B( x_{j_1 (n)}^0,  \ph( t_{j_1(n)}) + 2 t_n ^1 ) \subset B( x _{j_1 (n)}^0, t_{j_1(n)} )
$$
and using (\ref{4}) we get
$$
\begin{array}{l}
q_{j_1(n)} ^1 ( t_n ^1) \ges {\mu_{j_1(n)}}_{| B( x _{j_1(n)}^0, t _{j_1(n)})} \left( B( x_n ^1, t_n ^1) \right)
= \mu_{j_1(n)} \left( B( x_n ^1, t_n ^1) \right)
> q_{j_1(n)} \left( \ph( t_n^1) \right) - \frac{1}{ 2 j_1 (n)} .
\end{array}
$$
Since $ q_{j_1(n)} ^1 ( t_n ^1) \lra \al _1 < \al $ and $  q_{j_1(n)} \left( \ph( t_n^1) \right) \lra \al $ as	$n \lra \infty$
we see that 
the last inequality cannot hold for $n$ sufficiently large.  The claim (\ref{6}) is thus proven.

Replacing $ j_1$ by $ j_1( n_1 + \cdot)$, we may assume that (\ref{6}) holds for any $ n \in \N^*$.

\medskip

Let $q_n^2(t)$ be the concentration function of $ {\mu_{j_1(n)}}_{| B(x_n^1, \, t_n ^1)}.$
Let $ \tilde{q}^2(t) = \ds \limsup_{n \ra \infty} q_n^2 (t) .$
As previously, by Lemma \ref{e2}  there are a nondecreasing function $ q^2 : [0, \infty) \lra [0, M]$
and an increasing mapping  $ j_2 : \N^* \lra \N^*$ such that
\beq
\label{7}
q_{ j_2 (n)} ^2 (t) \lra q^2 (t) \qquad \mbox{ for all } t > 0 \mbox{ and }
\lim_ { t \ra \infty} q^2(t) > \lim_ { t \ra \infty} \tilde{q} ^2(t) - \frac{ \al }{2^3}.
\eeq

Let $ \al _2 = \ds \lim_ { t \ra \infty} q ^2(t).$
Since $ q_n ^2 \les q_{ j_1 (n)} $, it is clear that $ q^2 \les q $ and $ 0 \les \al _2 \les \al.$
If $ \al _2 = \al $ it follows from (\ref{5}) and Remark \ref{R2} (iv) that
$\left( {\mu _{ j_1(j_2(n))}}_ {| B( x_{j_2(n)} ^1, t_{j_2(n)} ^1 )}\right)_{n \ges 1} $
concentrates,
hence $\left( \mu_{j_1(j_2(n))} \right)_{n \ges 1}$,
$\left( x_{j_2(n)}^1 \right)_{n \ges 1}$ and   $\left( t_{j_2(n)} ^1 \right)_{n \ges 1}$  
satisfy the conclusion of Lemma \ref{L6} (relaxed).

If  $ \al _2 < \al$, we continue as above.
By Lemma \ref{e} (ii)-(iii)  there is an increasing  sequence $ t_n ^2 \lra \infty $
such that $ t_n ^2 < \frac 14 t_{j_2(n)} ^1  $ and $ q_{j_2(n)}^2 ( t_n ^2) \lra \al _2$.
Since $  t_n ^2 < \frac 14 t_{j_2(n)} ^1 < \frac{1}{4^2} t_{j_1 (j_2(n))}$, Lemma \ref{e} (iii) implies that
$$
q_{j_1( j_2(n))} ( t_n ^2) \lra \al  \qquad \mbox{ and } \qquad q_{j_1( j_2(n))} ( \ph(t_n ^2)) \lra \al \qquad
\mbox{ as } n \lra \infty.
$$
For each $n$ choose $ x_n ^2 \in \Xo$ such that
$$
\mu _{ j_1( j_2(n))} \left( B ( x_n ^2 , \ph(t_n ^2)) \right) > q_{j_1( j_2(n))} ( \ph(t_n ^2))  - \frac{1}{2^2 j_1( j_2(n))}.
$$
Obviously,
$$
q_{j_1( j_2(n))} ( t_n ^2) \ges \mu _{ j_1( j_2(n))} \left( B ( x_n ^2 , t_n ^2) \right)
\ges \mu _{ j_1( j_2(n))} \left( B ( x_n ^2 , \ph(t_n ^2)) \right)
>  q_{j_1( j_2(n))} ( \ph(t_n ^2))  - \frac{1}{2^2 j_1( j_2(n))},
$$
hence
$$
\mu _{ j_1( j_2(n))} \left( B ( x_n ^2 , t_n ^2) \right) \lra \al \qquad \mbox{ and } \qquad
\mu _{ j_1( j_2(n))} \left( B ( x_n ^2 , \ph( t_n ^2))	\right)  \lra \al \qquad \mbox{ as } n \lra \infty.
$$

There is $ n_2 \in \N^*$ such that for all $ n \ges n_2$,
\beq
\label{8}
B ( x_{j_1( j_2(n))}^0, \ph( t_{j_1( j_2(n))}))  \cap B( x_n ^2, t_n ^2) = \emptyset
\quad \mbox{ and } \quad
B ( x_{ j_2(n)}^1, \ph( t_{ j_2(n)} ^1 )) \cap B( x_n ^2, t_n ^2) = \emptyset.
\eeq
The proof of (\ref{8}) is similar to the proof of (\ref{6}) and we omit it.

Replacing $ j_2$ by $ j_2( \cdot + n_2) $ we may assume that (\ref{8}) holds for any $ n \in \N^*$.

\medskip

We continue the process inductively.
Suppose that we have constructed increasing mappings $ j_1, \dots,     j_k: \N^* \lra \N^*$ and for each 
$ i \in \{ 1, \dots, k \}$ we have an increasing sequence  $ (t_n^i)_{n \ges 1} $, 
a sequence of points
$ (x_n^i)_{n \ges 1}\subset \Xo$ and nondecreasing functions $ q_n^i, q^i: [0, \infty) \lra [0, M]$ 
with the following properties:

\medskip

{\bf (P1)}  $ t_n ^i \lra \infty $ as $ n \lra \infty, \;  $  $ t_n ^i < \frac 14 t_{j_i(n)} ^{ i-1} $
for all $ n \in \N^*$ and $ i \in \{ 1, \dots, k \}$, where $ t_n ^0 = t_n$. 

\medskip

As an obvious consequence of (P1) we have  
\beq
\label{9}
t_n^k < \frac 14 t_{j_k(n)} ^{k-1} < \frac{1}{4^2} t_{ j_{k-1}( j_k(n))} ^{k-2} < \dots 
< \frac{ 1}{ 4^{\ell}}  t_{j_{k-\ell +1}( \dots (j_k(n))\dots ) }^{ k - \ell } < \dots < \frac{1}{4^k} t _{ j_1( j_2(\dots ( j_k(n) \dots)))} .
\eeq

\medskip

{\bf (P2)}  $ q_n^i (t)$ is the concentration function of $ { \mu_{ j_1( j_2( \dots  (j_{i-1} (n)) \dots )) }}_{| B( x_n ^{i-1}, t_n^{i-1})}$,

\medskip

{\bf (P3)}  $ q_{j_i(n)} ^i( t) \lra q^i(t)$ as $ n \lra \infty$   and	   $ q^i(t) \lra \al _i $ as $ t \lra \infty$,
where 
\beq
 \ds \lim_ { t \ra \infty}   \left(   \limsup_{n \ra \infty} q_n ^i (t)  \right)    - \frac{ \al}{ 2^{i+1}} < \al _i < \al .
\eeq

\medskip

{\bf (P4)}  $ q_{j_i (n)}^i( t_n ^i) \lra \al _i$,    $ q_{j_i (n)}^i( \ph(t_n ^i) ) \lra \al _i$   	 as	  $n \lra \infty$ and
\beq \label{10}
\mu_{ j_1( \dots( j_i (n)) \dots )} \left( B( x_n ^i, \ph( t_n ^i))  \right) > q _{ j_1( \dots( j_i (n)) \dots )} (    \ph( t_n ^i))  - \frac{1}{ 2^i  j_1( \dots( j_i (n)) \dots )}.
\eeq

\medskip

{\bf (P5)}   We have  for all $ n \ges 1$, $\ell = 1, \dots, k $  and $  i = 1, \dots , \ell, $ 
\beq
\label{disjoint}
B \left( x_{j_i( j_{i+1} ( \dots ( j_{\ell} (n) \dots )) )} ^{i-1} ,
 \ph ( t_{j_i( j_{i+1} ( \dots ( j_{\ell} (n) \dots )) )} ^{i-1} ) \right)  \cap B( x_n ^{\ell}, t_n ^{\ell}) = \emptyset .
 \eeq

 Let $ q_n^{ k+1} (t)$ be the concentration function of $ {\mu _ {j_1 (\dots ( j_k (n)) \dots ) }}_{| B( x_n ^k, t_n ^k)} $. 
We have 
$$
\begin{array}{l}
q_{j_1( \dots ( j_k(n)) \dots )} ( t_{j_1( \dots ( j_k(n)) \dots )} ) \ges  q_{j_1( \dots ( j_k(n)) \dots )}  ( t_n ^k)  \qquad \qquad \mbox{ by (\ref{9}) }
\\
\\
\ges \mu_{j_1( \dots ( j_k(n)) \dots )}  ( B ( x_n ^k , t_n ^k)) \ges \mu_{j_1( \dots ( j_k(n)) \dots )}  ( B ( x_n ^k , \ph( t_n ^k)))
\\
\\
\ges q_{j_1( \dots ( j_k(n)) \dots )}  ( \ph( t_n^k)) - \frac{1}{2^k j_1( \dots ( j_k(n)) \dots )}   \qquad  \qquad \mbox{ by (\ref{10}) }.
\end{array}
$$
Hence using (\ref{conv}), the fact that $ q_n ( t_n ) \lra \al $ and Lemma \ref{e} (iii) we get 
$$ 
\mu_{j_1( \dots ( j_k(n)) \dots )}  \left( B ( x_n ^k ,  t_n ^k) \right)  \lra \al \quad \mbox{ and } \quad 
\mu_{j_1( \dots ( j_k(n)) \dots )}  \left( B ( x_n ^k , \ph( t_n ^k)) \right) \lra \al \quad \mbox{  as } n \lra \infty.
$$ 
It is clear that for all $ t \ges t_n^k$ we have 
$ q_n^{ k+1} (t) =  \mu_{j_1( \dots ( j_k(n)) \dots )}  ( B ( x_n ^k ,  t_n ^k))$
and the last quantity tends to $ \al $  as $ n \lra \infty$. 

Let $ \tilde{q}^{ k+1}(t) = \ds \limsup_{n \ra \infty} q_n^{k+1}(t)$.
As above, using Lemma \ref{e2}  we find  an increasing mapping  $ j_{k+1}  : \N^* \lra \N^*$ and 
a nondecreasing function $ q^{k+1} : [0, \infty) \lra [0, M]$ such that
\beq
\label{12}
q_{ j_{k+1} (n)} ^{k+1} (t) \lra q^{k+1} (t) \qquad \mbox{ for all  } t > 0 \mbox{ and }
\lim_ { t \ra \infty} q^{k+1}(t) > \lim_ { t \ra \infty} \tilde{q} ^{k+1} (t) - \frac{ \al }{2^{k+2}}.
\eeq
Let  $ \al _{k +1} = \ds \lim _{t \ra \infty} {q} ^{k+1} (t)$. 
From the above it is clear that $ \al_{k +1} \in [0, \al ]$. 
If $ \al_{k +1} = \al $ the conclusion of Lemma  \ref{L6} (relaxed) holds for the subsequences
$\left( \mu_{ j_1( \dots j_{k+1}(n)\dots )} \right)_{n \ges 1}$, 
$\left( x_{j_{k+1}(n) }^k \right)_{n \ges 1} $ and $\left( t_{j_{k+1}(n) }^k \right)_{n \ges 1} $.

If $ \al _{k +1} < \al$ we perform another step in our extraction process. 
By Lemma \ref{e} (ii)-(iii)  there is an increasing  sequence $ t_n ^{k+1} \lra \infty $  such that 
\beq
\label{13}
 t_n ^{k+1} < \frac 14 t_{j_{k+1}(n)} ^k  \qquad  \mbox{  and } \qquad q_{j_{k+1}(n)}^{k+1} ( t_n ^{k+1}) \lra \al _{k+1} \qquad \mbox{ as } n \lra \infty.
\eeq
Since $  t_n ^{k+1} < \frac 14 t_{j_{k+1}(n)} ^k < \frac{1}{4^{k+1}} t_{j_1( \dots j_{k+1}(n)) \dots )}$, 
Lemma \ref{e} (iii) implies that
$$
q_{j_1 ( \dots ( j_{k+1}(n)) \dots )} ( t_n ^{k+1}) \lra \al  \qquad \mbox{ and } \qquad q_{j_1(\dots ( j_{k+1} (n)) \dots )} ( \ph(t_n ^{k+1})) \lra \al \qquad
\mbox{ as } n \lra \infty.
$$
For each $n$ we choose $ x_n ^{k+1} \in \Xo$ such that
$$
\mu _{ j_1(\dots ( j_{k+1}(n)) \dots) } \left( B ( x_n ^{k+1} , \ph(t_n ^{k+1})) \right) 
> q_{j_1( \dots (  j_{k+1}(n)) \dots ) } ( \ph(t_n ^ {k+1}))  - \frac{1}{2^{k+2} j_1( \dots (j_{k+1} (n)) \dots )}.
$$
Since 
$$
q_{j_1 ( \dots ( j_{k+1}(n)) \dots )} ( t_n ^{k+1})  
\ges 
\mu _{ j_1(\dots ( j_{k+1}(n)) \dots) } \left( B ( x_n ^{k+1} , t_n ^{k+1}) \right) 
\ges 
\mu _{ j_1(\dots ( j_{k+1}(n)) \dots) } \left( B ( x_n ^{k+1} , \ph(t_n ^{k+1})) \right), 
$$
we have 
\beq
\label{13bis}
\begin{array}{l}
\mu _{ j_1(\dots ( j_{k+1}(n)) \dots) } \left( B ( x_n ^{k+1} , t_n ^{k+1}) \right) \lra \al  \qquad \mbox{ and } 
\\
\\
\mu _{ j_1(\dots ( j_{k+1}(n)) \dots) } \left( B ( x_n ^{k+1} , \ph(t_n ^{k+1})) \right)  \lra \al 
\qquad \mbox{ as } n \lra \infty.
\end{array}
\eeq
We show that there is $ n_{k +1} \in \N$  such that for all $ n \ges n_{k+1}, $
\beq
\label{disjoints}
\forall i \in \{ 1, \dots, k+1 \}, \qquad 
B \left(x_{j_i( \dots (j_{k+1} (n)) \dots ) } ^{i-1} , \ph ( t _{j  _i( \dots (j_{k+1} (n)) \dots ) } ^{i-1} ) \right)  
\cap B( x_n ^{k+1}, t_n^{k+1} ) = \emptyset. 
\eeq
Indeed, suppose that for some $ i \in \{ 1, \dots, k+1 \} $ the intersection in (\ref{disjoints}) is not empty. 
Then 
\beq
\label{15}
\begin{array}{lcl}
B( x_n ^{k+1}, t_n^{k+1} ) 
& \subset & 
B \left(x_{j_i( \dots (j_{k+1} (n)) \dots ) } ^{i-1} , \ph ( t _{j_i( \dots (j_{k+1} (n)) \dots ) } ^{i-1} )  + 2 t_n^{k+1} \right)  
\\
 & \subset &  
B \left(x_{j_i( \dots (j_{k+1} (n)) \dots ) } ^{i-1} ,  t _{ j _i( \dots (j_{k+1} (n)) \dots ) } ^{i-1} ) \right)  
\end{array}
\eeq
because $ \ph (s) \les \frac s2$ and $ t_n^{k+1} < \frac 14 t _{ j _i( \dots (j_{k+1} (n)) \dots ) } ^{i-1}  .$
If $ i \les k $ we get 
$$
\begin{array}{l}
q_{j_i( j_{i+1} ( \dots ( j_{k+1} (n)) \dots  ) ) }^ i  \left( t_{j_{i+1} ( \dots (( j_{k+1} (n)) \dots ) ) } ^ i \right) 
\\
\\
\ges q_{j_i( j_{i+1} ( \dots ( j_{k+1} (n)) \dots  ) ) }^ i  ( t_n ^{ k+1 } ) 
\qquad \qquad \mbox{ by (\ref{9}) and (\ref{13})}
\\
\\
\ges { \mu _{j_1 ( \dots ( j_{k+1} (n)) \dots ) } } _{ \Big| B\left( x_ {j_i( j_{i+1} ( \dots ( j_{k+1} (n)) \dots  ) ) }^ {i-1} ,  t_{j_i( j_{i+1} ( \dots ( j_{k+1} (n)) \dots  ) ) }^ {i-1}  \right)}  
\left( B( x_n ^{k +1}, t _n ^{k +1}  )\right) 
\\
\\ 
= \mu _{j_1 ( \dots ( j_{k+1} (n)) \dots ) }  \left( B( x_n ^{k +1}, t _n ^{k +1}  )\right)  \qquad \qquad \mbox{ by (\ref{15})}.
\end{array}
$$
As $ n \lra \infty$ 
the right-hand side of the above inequality tends to $ \al $ by (\ref{13bis}) and the left-hand side tends to $ \al _i < \al $ 
by (P4), 
hence the inequality  may be true  only for finitely many $n$'s. 

If $ i = k +1$ we have
$B( x_n  ^{k+1}, t_n  ^{ k +1 } ) \subset B \left( x_{j_{k+1}(n)} ^k, t _{j_{k+1}(n)} ^k \right), $ hence
$$
\begin{array}{l}
q_{j_{k+1} (n)} ^{k+1} ( t_n^{k+1}) \ges {\mu_{j_1(  \dots ( j_{k+1} (n )) \dots ) } } _{\big|   B \left( x_{j_{k+1}(n)} ^k, t _{j_{k+1}(n)} ^k \right) }
\left( B( x_n  ^{k+1}, t_n  ^{ k +1 } )   \right) 
\\
= \mu_{j_1(  \dots ( j_{k+1} (n )) \dots ) }  \left( B( x_n  ^{k+1}, t_n  ^{ k +1 } )   \right) .
\end{array}
$$
Since $ q_{j_{k+1} (n)} ^{k+1} ( t_n^{k+1})  \lra \al _{k +1} < \al $ by (\ref{13}) and 
$ \mu_{j_1(  \dots ( j_{k+1} (n )) \dots ) }  \left( B( x_n  ^{k+1}, t_n  ^{ k +1 } )   \right)  \lra \al $ as $ n \lra \infty$ by (\ref{13bis}), 
the above inequality cannot hold for  $n$ sufficiently large. 

We conclude that  there is $ n_{k+1} \in \N$ such that (\ref{disjoints}) holds for all $n \ges n_{k +1}$. 

Replacing $ j_{k+1}$ by $j_{k+1}( \cdot + n_{k+1}) $ we may assume that (\ref{disjoints}) holds for all $n  \in \N^*$. 

The mappings $ j_1, \dots, j_{k +1}$, the functions $ q_n^i $ and  $ q^i$ and the sequences $ (x_n^i)_{n \ges 1}$, $ (t_n^i)_{n \ges 1}$ 
with $ i \in \{ 1, \dots ,  k+1 \}$   satisfy the properties (P1) - (P5) above and this finishes our induction.
We conclude that  at step $k$ either 
we find a subsequence of $ (\mu_n)_{n \ges 1}$ and a sequence of balls satisfying the conclusion of Lemma \ref{L6} (relaxed), or
we are able to complete the step $ k+1$. 

\medskip

We claim that the above process has to stop after a finite number of steps. 
Indeed, assume that we have completed $k$ steps. It follows immediately from (\ref{disjoint}) that for each $ n \in \N^*$ the balls 
$$
B \left( x_{j_1( \dots ( j_k( n)) \dots ) } ^0, \ph ( t_ {j_1( \dots ( j_k( n)) \dots ) } ) \right),  \dots , 
B \left( x_{j_{\ell +1} ( \dots ( j_k( n)) \dots ) } ^{\ell}, \ph ( t_ {j_{\ell +1} ( \dots ( j_k( n)) \dots ) } ^{\ell} ) \right), 
\dots , 
B( x_n^k, t_n ^k) 
$$
are disjoint, where $ \ell = 1, \dots , k-1$. 
From the choice of $ x_n ^0$, (\ref{4}) and (\ref{10}) we have 
$$
\begin{array}{l}
\mu_{j_1 ( \dots ( j_k (n)) \dots ) }  \left( 
B \left( x_{j_{\ell +1} ( \dots ( j_k( n)) \dots ) } ^{\ell}, \ph ( t_ {j_{\ell +1} ( \dots ( j_k( n)) \dots ) } ^{\ell} ) \right) \right) 
\\
\qquad 
>   q_{j_1( \dots ( j_k(n))) \dots )} \left(  \ph ( t_ {j_{\ell +1} ( \dots ( j_k( n)) \dots ) } ^{\ell} ) \right) 
- \frac{1}{2^{\ell} j_1( \dots ( j_k(n))) \dots )} \quad \mbox{ for } \ell = 0, 1, \dots, k-1, \; \;  \mbox{ and } 
\end{array}
$$
$$
\mu_{j_1 ( \dots ( j_k (n)) \dots ) }  \left(  B( x_n ^k, \ph (t_n ^k)) \right) 
> q_{j_1( \dots ( j_k(n))) \dots )} \left(  \ph ( t_n ^k) \right) - \frac{1}{2^{k } j_1( \dots ( j_k(n))) \dots )} .
$$
Since the balls are disjoint, summing up the above inequalities we get 
$$ 
\begin{array}{l}
\mu_{j_1 ( \dots ( j_k (n)) \dots ) }   (\Xo) 
\\
\ges 
\ds \sum_{\ell = 0} ^{ k-1}  q_{j_1( \dots ( j_k(n))) \dots )} \left(  \ph ( t_ {j_{\ell +1}( \dots ( j_k( n)) \dots ) } ^{\ell} ) \right)  
+ q_{j_1( \dots ( j_k(n))) \dots )} \left(  \ph ( t_n ^k) \right)
 -\frac{2}{ j_1( \dots ( j_k(n))) \dots )} . 
\end{array}
$$
We rememeber that $  \ph ( t_ {j_{\ell +1}( \dots ( j_k( n)) \dots ) } ^{\ell} )  \lra \infty $ as $ n \lra \infty $ for each $ \ell$, 
then we take the  $\ds \limsup$ in the above inequality and use Lemma \ref{e} (i) to get  
$
M \ges ( k+ 1) \al .
$
This implies that $k$ has to remain  bounded, 
which means that  we may perform only a finite number of steps in the previous extraction process. 
The relaxed version of Lemma \ref{6} is thus proven. 

\medskip

It remains to show  that the subsequences may be choosen in such a way that 
$( \mu_{n_k})_{k \ges 1}$ concentrates around $ (x_k )_{k \ges 1}$. 

Let $ { \ph} (t)  = \frac{t}{12}$. 
There exist a subsequence $( \mu_{n_k})_{k \ges 1}$, 
a sequence $ (z_k)_{k \ges 1} \subset \Xo$ and an increasing sequence $ \tilde{R}_k \lra \infty $ satisfying 
the relaxed version of Lemma \ref{L6}  with 
$ z_k,  \tilde{R}_k$ instead of $ x_k$ and  $R_k$, respectively. In particular,   the sequence of measures 
$\left({\mu_{n_k}}_{ | B( z_k, \tilde{R}_k)} \right)_{k \ges 1}$ 
concentrates around some sequence of points $ (x_k)_{k \ges 1}$. 

Using Remark \ref{R2} (i) we get 
$$
\begin{array}{l}
\mu_{n_k} \left( B( z_k, \tilde{R}_k ) \setminus B\left( x_k, \frac{ \tilde{R}_k}{12}\right) \right) 
= {\mu_{n_k}}_{ \big| B( z_k, \tilde{R}_k)} \left( \Xo \setminus B\left( x_k, \frac{ \tilde{R}_k}{12} \right) \right)   \lra 0 
\qquad \mbox{ as } k \lra \infty, 
\end{array}
$$
and this implies
$
\mu_{n_k} \left( B\left( z_k, \tilde{R}_k \right) \cap B\left( x_k, \frac{ \tilde{R}_k}{12} \right) \right) \lra \al .
$
On the other hand we have
$$
\begin{array}{l}
\mu_{n_k} \left( B \left( z_k, \frac{\tilde{R}_k}{12} \right) \right) \ges
\mu_{n_k} \left( B \left( z_k, \frac{ \tilde{R}_k}{12} \right)  \cap B\left( x_k, \frac{ \tilde{R}_k}{12} \right) \right)
\\
\\
\ges
\mu_{n_k} \left( B \left( z_k, \tilde{R}_k \right) \cap B\left( x_k, \frac{ \tilde{R}_k}{12}  \right) \right)  - 
\mu_{n_k} \left( B\left( z_k, \tilde{R}_k\right) \setminus  B \left( z_k, \frac{ \tilde{R}_k}{12} \right) \right) .
\end{array}
$$
Using (\ref{annulus})  we infer  that 
$ \mu_{n_k} \left( B \left( z_k, \frac{\tilde{R}_k}{12} \right)  \cap B \left( x_k, \frac{ \tilde{R}_k}{12} \right) \right) \lra \al .$
In particular, for all $k$ sufficiently large we have 
$ B\left( x_k, \frac{ \tilde{R}_k}{12}  \right) \cap B \left( z_k, \frac{ \tilde{R}_k}{12} \right) \neq \emptyset , $ 
hence $ d( x_k, z_k) \les \frac{ \tilde{R}_k}{6}.$
Let $ R_k = \frac 12 \tilde{R}_k. $  Then 
$$
\begin{array}{l}
B\left( z_k, \frac{ \tilde{R}_k }{12}  \right) \subset B\left( x_k,  \frac{ \tilde{R}_k}{4}  \right)
=  B\left( x_k,  \frac{{R}_k}{2} \right) \subset  B( x_k,  {R}_k )
\subset B( z_k, \tilde{R}_k). 
\end{array}
$$
Hence $ \mu_{n _k} ( B( x_k,  {R}_k)   ) \lra \al$ as $ k \lra \infty$, 
the sequence of measures
$\left(  {\mu _{n _k}} _{ |  B( x_k,  {R}_k)  } \right)_{k \ges 1} $ 
concentrates around $ (x_k)_{k \ges 1} $ (cf. Remark \ref{R2} (i))  and  Lemma \ref{L6} is proven. 
\hfill
$\Box$

\section{Proof of Theorem \ref{main} }
\label{Prmain}

Let $ q_n(t) $ be the concentration function of $ \mu_n$ and let 
$$
\al _0:= \lim_{t \ra \infty} \left( \limsup_{ n \ra \infty} q_n (t) \right) .
$$
If $ \al _0= 0 $ then $ (\mu _n)_{ n \ges 1}$ is a vanishing sequence and the conclusion of  Theorem \ref{main}  holds. 

\medskip

From now on we assume that $ \al _0 > 0 $ and we use the extraction procedure described below. 

\medskip

{\it Step 1. } Let $ (\nu _n)_{ n \ges 1}$ be a sequence of positive Borel measures on $\Xo$ 
with concentration functions $ q_n $ such that $\nu _n (\Xo)$ is bounded and 
\beq
\label{alf}
\lim_{t \ra \infty } \left( \limsup_{ n \ra \infty } q_n (t) \right) = \al  >0. 
\eeq
Fix $ s > 0$, $ \e _1, \e _2 > 0$ and an increasing sequence of positive numbers $ (t_n)_{n \ges 1}$,  $ t_n \lra \infty.$

We show that there exist an increasing mapping $ j : \N^* \lra \N^*$, a sequence of points $(x_n )_{n \ges 1}\subset \Xo$ 
and an increasing sequence $ (R_n)_{n \ges 1}$, $ R_n \lra \infty $ with  the following properties: 

\medskip

(p1) $ \quad $ $ q_{j (n)} ( s) < 2 \al $ for any $n$, 

\medskip

(p2) $ \quad $  $R_n \les t_{j(n)} $  for all $n$, 

\medskip

(p3) $ \quad $  $ 0 < \nu _{ j(n)} \left( B( x_n, \ph( R_n)) \right) \lra \al '$ as $ n \lra \infty$, where $ \al ' > \al - \e_1$, 

\medskip

(p4)  $ \quad $  $ \nu _{ j(n)} \left( B( x_n, 2 R_n) \setminus  B( x_n, \ph( R_n)) \right) < \frac{ \e _2}{2^n }$, 

\medskip

(p5) $ \quad $  $ \left( {\nu _{ j(n)}}_{|  B( x_n, 2 R_n) } \right) $ concentrates around $( x_n)_{n \ges 1}. $

\medskip

By (\ref{alf}) we have $\ds \limsup_{ n \ra \infty } q_n (t) \les \al $ for all $ t >0$. Eliminating a finite number of terms we may assume that 
$ q_n ( s ) <2 \al $ for all $ n$. 

Using Lemma \ref{e2} we find an increasing mapping $ \kappa_1 : \N^* \lra \N^*$, a nondecreasing function
 $ q : [0, \infty) \lra [0, \infty) $  such that 
$ \ds \lim_{n \ra \infty } q_{\kappa_1 (n)  }(t) = q(t)$  for any $ t > 0$ and 
 $ \ds \lim_{t \ra \infty} q(t) > \al - \e _1$.

Let $ \al '= \ds \lim_{t \ra \infty} q(t) .$
We may assume that $ \ds \lim_{ n \ra \infty} q_{\kappa_1 (n)  } ( t_{\kappa_1 (n) } ) = \al '$ 
(otherwise we replace $( t_{\kappa_1 (n) } )_{n \ges 1}$ by a  sequence $( \tilde{t}_{\kappa_1 (n) } )_{n \ges 1}$ 
which has this property and satisfies $  \tilde{t}_{\kappa_1 (n)} \les  t_{\kappa_1 (n)}  $; 
the existence of $( \tilde{t}_{\kappa_1 (n) } )_{n \ges 1}$  is guaranteed by Lemma \ref{e} (ii)-(iii)). 

Using Lemma \ref{L6} 
for the sequences
 $ \left(\nu_{ \kappa_1(n)} \right)_{n \ges 1}$ and $ \left(t_{ \kappa_1(n)} \right)_{n \ges 1}$
we find an increasing mapping $ \kappa_2 : \N^* \lra \N^*$, 
a sequence of points $ (y_n)_{ n \ges 1} \subset \Xo$ and a sequence $ R_n^{*} $ such that 
$  R_n^{*}  \les t_{\kappa_1 ( \kappa_2 (n))}$, $  R_n^{*}  \lra \infty$,  
$ \nu_{\kappa_1 ( \kappa_2 (n)) } \left( B\left( y_n, \ph \left( {  R_n^{*}}/{2} \right) \right) \right) \lra \al ', $
\beq
\label{16}
\beta_n := \nu_{\kappa_1 ( \kappa_2 (n)) } \left( B( y_n, R_n ^{*})  \setminus B\left( y_n, \ph \left( {  R_n^{*}}/{2} \right) \right) \right) \lra 0
\qquad \mbox{ as }  n \lra \infty
\eeq
and the sequence $ \left( { \nu_{\kappa_1 ( \kappa_2 (n)) }}_{\big|  B( y_n, R_n ^{*}) } \right)_{ n \ges 1}$ concentrates. 

We construct $ \kappa_3 : \N^* \lra \N^* $ inductively as follows: 
$ \kappa_3 (1) $ is the first integer satisfying $ \beta_{  k_1} <\frac 12$ (where $ \beta _n$ is defined in (\ref{16})). 
If $ \kappa_3 (1), \dots ,  \kappa_3 (n)$ have been constructed, we define $ \kappa_3 (n+1)$ as the first integer $ k_{n +1}$ verifying
$ k_{n + 1} > \kappa_3 (n)$, $ \beta _{ k_{n +1}} < \frac{1}{2 ^{ n+1}} $ and $ R_{ \kappa_3 (n) }^* < R_{k_{n+1}}^*.$

Let $ j = \kappa_1 \circ \kappa_2 \circ \kappa_3, $
$ \; x_n = y_{\kappa_3(n)}$ and $ R_n = \frac 12  R_{ \kappa_3 (n) }^*$. 
Then $j$, $( x_n)_{n \ges 1}$  and  $( R_n)_{n \ges 1}$  have  the properties (p1)$-$(p5) listed above. 

\medskip

We apply step 1 to $ (\mu_n )_{n \ges 1}$ with $ s = 1$, $ \e _1 = \frac{ \al _0}{2} $, $ \e _2 = 1$ and we get 
an increasing mapping $ j_1 : \N^* \lra \N^*$, a sequence of points $ (z_n ^1)_{n \ges 1} \subset \Xo$ and 
an increasing sequence $ (R_n ^1)_{n \ges 1} $  with the properties (p1)$-$(p5) above. 
Denote 
$$
 m_1 = {\ds \lim_{n \ra \infty} }\mu_{ j_1 (n)} ( B( z_n ^1, \ph( R_n ^1))) \in \Big( \frac{\al_0}{2}, \al _0 \Big].
$$

Let $ \mu_n ^1 = {\mu_{j_1(n) }}_{\big| \Xo \setminus B( z_n ^1, R_n ^1) } $ and let $ q_n ^1$ be the concentration function of 
$  \mu_n ^1$. 
Denote 
\beq
\label{17}
\al _1 = \lim_{t \ra \infty} \left( \limsup_{ n \ra \infty} q_n ^1 (t) \right) .
\eeq
Since $ q_n ^1 ( t) \les q_{j_1(n)}(t)$ it is clear that $ 0 \les \al _1 \les \al _0$. 
If $ \al _1  = 0 $ conclusion (i) in Theorem \ref{main} is satisfied (with $ k =1$) by $j_1$ and the balls $B( z_n ^1, R_n ^1)$. 
If $ \al _1 > 0 $ we proceed to step 2. 

\medskip 

{\it Step 2. } We apply step 1 to the sequence $ (\mu_n ^1)_{ n \ges 1}$ with $ s = 2$, 
$ \e _1 = \frac{ \al _1}{2}$,  $ \e _2 = \frac 12$  and $ t_n = \frac 12 R_n ^1$. 
We obtain an increasing mapping $ j_2 : \N^* \lra \N^*$, a sequence if points $ ( z_n ^2)_{n \ges 1} \subset \Xo$, 
and an increasing sequence $ (R_n^2)_{n \ges 1} \subset (0, \infty)$    such that 

\medskip 

$\bullet$ $\quad $ $ R_n ^2 \lra \infty $ and $ R_n ^2 \les \frac 12 R_{j_2 (n)} ^1 $, 

\medskip 

$\bullet$ $\quad $  $ \mu_{ j_2(n)} ^1 \left( B ( z_n ^2, \ph  (R_n^2 ) ) \right) \lra  m_2 $, where $ \frac{ \al _1}{2} < m_2 \les \al _1$, 

\medskip 

$\bullet$ $\quad $ 
$  \mu_{ j_2(n)} ^1 \left( B ( z_n ^2, 2 R_n ^2 )  \setminus  B ( z_n ^2, \ph  (R_n ^2) ) \right)  < \frac{1}{ 2^{n+1}}. $

\medskip 

We show that
\beq
\label{18} 
B( z_n ^2, R_n ^2) \cap B( z_{ j_2(n)} ^1, R_{j_2(n)}^1 ) = \emptyset 
\eeq
for all $n$ sufficiently large.
Indeed, if the two balls intersect we  have
$$
B( z_n ^2, R_n ^2)  \subset B( z_{ j_2(n)} ^1, R_{j_2(n)}^1 + 2 R_n^2) 
\subset B( z_{ j_2(n)} ^1, 2 R_{j_2(n)}^1 ) , 
$$
hence 
$$
\mu_{j_2(n)} ^1 ( B( z_n ^2 , R_n^2 ))  \les \mu_{j_2(n)} ^1 \left(  B( z_{ j_2(n)} ^1, 2 R_{j_2(n)}^1 )  \right)
= \mu _{ j_1( j_2(n)) } \left( B( z_{ j_2(n)} ^1,2 R_{j_2(n)}^1 ) \setminus B( z_{ j_2(n)} ^1, R_{j_2(n)}^1 ) \right).
$$
The left hand side in the above inequality tends to $  m_2>0$ and the right hand side tends to $0$ as $ n \lra \infty$, 
consequently there is $ n_1 \in \N$ such that the inequality cannot hold for  $ n \ges n_1$.
Therefore  (\ref{18}) is true for $ n \ges n_1$. 
Replacing $ j_2 $ by $ j_2 (\cdot + n_1) $ we assume from now on that (\ref{18}) holds for all $ n \in \N^*$. 

Let 
$$
\mu_n ^2 = { \mu_{ j_2(n)} ^1 }_{ \big| \Xo \setminus B( z_n ^2,  R_n ^2) } 
= {\mu_{ j_1( j_2(n))} }_{\big| \Xo \setminus \left( B( z_{ j_2(n)}^1,  R_{j_2(n)}^1) \cup B( z_n ^2,  R_n ^2) \right) }. 
$$
Let $ q_n ^2 $ be the concentration function of $ \mu_n ^2$ and let 
$$
\al _2 = \lim_{t \ra \infty} \left( \limsup _{n \ra \infty} q_n ^2( t) \right) . 
$$
Since $ q_n ^2 \les q_{j_2(n)} ^1 $ we have obviously $ 0 \les \al _2 \les \al _1$. 
If $ \al _2 = 0 $ conclusion (i) in Theorem \ref{main} is satisfied with $ k = 2$, $ j = j_1 \circ j_2$, 
$ x_n ^1 = z_{ j_2(n)} ^1$, $ x_n ^2 = z _n ^2$, $ r_n ^1 = R_{j_2(n)}^1$, $ r_n ^2 = R_n ^2$. 

If $ \al _2 > 0 $ we proceed to step 3, which consists in applying step 1 to the sequence $ ( \mu_n^2)_{n \ges 1}$ with $ s =3$, 
$ \e _1 = \frac{ \al _2}{2}$, $ \e _2 = \frac{1}{2 ^2}, $ $ t_n = \frac 12 R_n ^2$, and so on. 

\medskip

We continue the above process inductively. 
Assume that we have completed $k $ steps and we have found $k$ increasing mappings $ j_1, \dots, j_k : \N^* \lra \N^*$, 
increasing sequences of positive numbers $ (R_n^1)_{n \ges 1}, \dots  (R_n^k)_{n \ges 1}$ that tend to infinity, 
 sequences of points $ ( z_n ^1)_{n \ges 1}, \dots ,  ( z_n ^k)_{n \ges 1} \subset \Xo$,
$k$ sequences of  measures $ (\mu_n ^i )_{n \ges 1}$ with concentration functions $ q_n ^i$, $ i \in \{ 0, \dots, k -1 \}$ 
(where $ \mu_n ^0 = \mu_n$, $ q_n ^0 = q_n$) 
and positive numbers $  \al _0 \ges \al _1 \ges \dots \ges \al _{ k-1} $, $ m_1, \dots, m_k $ 
satisfying  the following properties: 

\medskip 

\noindent
(H1) $\quad $ 
$ R_n ^k \les \frac 12 R_{ j_k (n)} ^{ k-1} \les \frac{1}{2^2} R_ { j_{ k-1} ( j_k (n)) } ^{ k-2} \les \dots \les
\frac{1}{ 2^{ k-1} } R_{ j_2( \dots ( j_k(n)) \dots ) } ^1. $

\medskip 

\noindent
(H2) $\quad $ 
$ \begin{array}{rcl}
\mu_n ^{\ell }&  = & {\mu_{ j_{\ell}(n) }^{\ell -1} }_{\big| \Xo \setminus B(  z_n ^{\ell}, R_n^{\ell}) }
\\
 & = & { \mu _{ j_1( \dots ( j_{\ell }(n )) \dots ) } } _{ \big| \Xo \setminus 
\left( B(  z_n ^{\ell}, R_n^{\ell}) \cup B\left( z_{ j_{\ell} (n) }^{ \ell -1} , R_{ j_{\ell} (n) }^{ \ell -1}\right) \cup \dots \cup 
B\left( z_{j_2( \dots ( j_{\ell }(n)) \dots )} ^1 , R _{j_2( \dots ( j_{\ell }(n)) \dots )} ^1 \right) \right) } 
\end{array} $

\medskip

\noindent
(H3) $\quad  $  
$ \ds \lim_{t \ra \infty} \left( \lim_{n \ra \infty} q_n ^{\ell -1} (t ) \right) \! = \! \al_{\ell -1} $ and 
$ \mu_{ j_{\ell}(n) } ^{ \ell -1 } \left( B( z_n ^{\ell}, \ph ( R _n ^{\ell})  ) \right) \lra m _{\ell }$, where 
$ \frac{ \al _{\ell -1}}{2} \! < \!  m _{\ell } \! \les \! \al_{\ell -1}$. 

\medskip

\noindent
(H4) $ \quad  $
$ \mu_{ j_{\ell}(n) } ^{ \ell -1 } \left( B( z_n ^{\ell}, 2 R _n ^{\ell}   )  \setminus B( z_n ^{\ell}, \ph ( R _n ^{\ell})  ) \right) 
\les \frac{ 1}{ 2^{ n + \ell -1} } $ for all $n$ and $ \ell = 1, \dots,  k$, where $ \mu_n ^0 = \mu _n$. 

\medskip 

\noindent
(H5)  $\quad $  The balls 
$$
B( z_n ^k, R_n ^k), \; 
B \left( z_{ j_k (n)} ^{ k-1}, R _{ j_k (n)} ^{ k-1} \right) \! , \; 
B \left( z_{ j_{ k-1} (j_k (n)) } ^{ k-2}, R _{ j_{k-1} (j_k (n))} ^{ k-2} \right)\! , \dots ,
B\left( z_{j_2( \dots ( j_{k }(n)) \dots )} ^1 , R _{j_2( \dots ( j_{k }(n)) \dots )} ^1 \right)  
$$
are all disjoint. 

\medskip 

\noindent
(H6)  $\quad $ The sequence 
$\left( {\mu_ { j_1( \dots ( j_{\ell}(n))\dots ) }}_{\big| B( z_n ^{\ell}, 2 R_n ^{\ell} ) } \right)_{n \ges 1}$ 
concentrates around $ (z_n^{\ell} )_{n \ges 1}.$

\medskip 

\noindent
(H7)  $\quad $   
$ m_1 + \dots + m_k   \les M. $

\medskip 

\noindent
(H8)  $\quad $ $ q_{ j_{\ell}(n)} ^{ \ell -1 } (\ell ) \les 2 \al _ {\ell -1}  $ for all $n$ and $ \ell = 1, \dots, k$.

\medskip

Let $ \mu_n ^k = {\mu_{ j_{k}(n) }^{k -1} }_{\big| \Xo \setminus B(  z_n ^{k}, R_n^{k}) } .$
Denote by $ q_n ^k $ the concentration function of $ \mu_n ^k  $ and  let
$
\al _k = \ds \lim_{t \ra \infty} \left( \limsup _{n \ra \infty} q_n ^k( t) \right) . 
$
Property (H2) above implies that $ q_n ^k \les q_{ j_{k} (n)} ^{ k-1}$, hence $ 0 \les \al _k \les \al _{k-1}.$
If $ \al _k = 0 $
conclusion (i) in Theorem \ref{main} is satisfied by $ j = j_1 \circ \dots \circ j_k$ and the balls in  (H5) above. 

\medskip

If $ \al _ k > 0 $ we proceed to step $ k+1$ which consists in applying step 1 to the sequence of measures 
$ (\mu_n ^k)_{ n \ges 1}$ with $ s = k+1$, 
$ \e _1 = \frac{ \al _k}{2}$, $ \e _2 = \frac{1}{ 2^k}$ and $ t_n = \frac 12 R_n^k$. 
We find an increasing mapping $ j_{ k + 1} : \N^* \lra \N^*$, a sequence of points $ (z_{n}^{k+1})_{ n \ges 1}$
 and an increasing  sequence of positive numbers $ (R_n^{ k+1} )_{n \ges 1}$ with the properties (p1)$-$(p5) in step 1. 
It is clear that (H1)-(H4), (H6) and (H8) hold with $ k+1$ instead of $k$.  
We claim that there is $ n_k \in \N$ such that for any  $ n \ges n_k$ and any $ \ell \in \{ 1, \dots, k \}$   we have 
\beq 
\label{19}
B( z_n ^{ k+1} , R_n^{ k+1}) \cap 
B\left( z_{j_{\ell +1} ( \dots ( j_{k +1}(n)) \dots )} ^{\ell}  , R _{j_{\ell +1} ( \dots ( j_{k +1}(n)) \dots )} ^{\ell} \right)
= \emptyset. 
\eeq
Indeed, if the intersection is not empty  using (H1) we get 
\beq
\label{20}
\begin{array}{rcl}
B( z_n ^{ k+1} , R_n^{ k+1}) & \subset & 
B\left( z_{j_{\ell +1} ( \dots ( j_{k +1}(n)) \dots )} ^{\ell}  , R _{j_{\ell +1} ( \dots ( j_{k +1}(n)) \dots )} ^{\ell} + 2 R_n^{ k+1}) \right)
\\
& \subset &  B\left( z_{j_{\ell +1} ( \dots ( j_{k +1}(n)) \dots )} ^{\ell}  , 2 R _{j_{\ell +1} ( \dots ( j_{k +1}(n)) \dots )} ^{\ell} ) \right).
\end{array}
\eeq
From properties (p3) and (p4) in step 1 it follows that 
$$
\mu_{ j_{ k+1} (n)} ^k \left( B( z_n ^{ k+1} , R_n^{ k+1}) \right) 
\lra m_{k+1}  \in \Big( \frac{ \al _k}{2}, \al _k \Big] 
\qquad \mbox{ as } n \lra \infty. 
$$
On the other hand using (H2) and (H4) we find 
$$
\begin{array}{l}
\mu_{ j_{ k+1} (n)} ^k  \left(  B\left( z_{j_{\ell +1} ( \dots ( j_{k +1}(n)) \dots )} ^{\ell}  , 2 R _{j_{\ell +1} ( \dots ( j_{k +1}(n)) \dots )} ^{\ell} ) \right) \right)
\\
\\
\les \mu_{ j_{\ell +1} ( \dots ( j_{ k+1} (n)) \dots )} ^{ \ell} 
\left(  B\left( z_{j_{\ell +1} ( \dots ( j_{k +1}(n)) \dots )} ^{\ell}  , 2 R _{j_{\ell +1} ( \dots ( j_{k +1}(n)) \dots )} ^{\ell} ) \right) \right)
\\
\\
= \mu_{ j_{\ell } ( \dots ( j_{ k+1} (n)) \dots ) }^{ \ell -1} 
\left(  B\left( z_{j_{\ell +1} ( \dots ( j_{k +1}(n)) \dots )} ^{\ell}  , 2 R _{j_{\ell +1} ( \dots ( j_{k +1}(n)) \dots )} ^{\ell} ) \right) 
\right.
\\
\left.
\qquad \qquad \qquad \qquad \qquad 
\setminus B\left( z_{j_{\ell +1} ( \dots ( j_{k +1}(n)) \dots )} ^{\ell}  ,  R _{j_{\ell +1} ( \dots ( j_{k +1}(n)) \dots )} ^{\ell} ) \right)\right)
\\
\\
\les \frac{1}{2^{j_{\ell +1} ( \dots ( j_{ k+1} (n)) \dots ) + \ell -1}} \lra 0 \qquad \mbox{ as } n \lra \infty. 
\end{array}
$$
We conclude that the  (\ref{20}) cannot be true if  $n$ is sufficiently large and  the claim (\ref{19}) is  proven. 
Replacing $ j_{ k +1} $ by $  j_{ k +1} ( \cdot + n_{ k } )$ we may assume that (\ref{19}) holds for all $ n \in \N^*$. 
In particular, we see that (H5) is satisfied at level $k+1$. 

Let $ \mu _n ^{ k + 1} = {\mu_{ j_{k+1 }(n) }^{k} }_{\big| \Xo \setminus B(  z_n ^{k+1}, R_n^{k+1}) } $. 
For $ \ell = 1, \dots , k+1 $ we have 
$$
\begin{array}{l}
\mu_{ j_1( \dots ( j_{k+1} ( n)) \dots )} \left( B \left( z_{j_{\ell +1}  ( \dots ( j_{k+1} ( n)) \dots )}^{\ell},  R _{j_{\ell +1}  ( \dots ( j_{k+1} ( n)) \dots )}^{\ell} \right) \right)
\\ 
\ges
\mu _{j_{\ell } ( \dots ( j_{k+1} ( n)) \dots )} ^{ \ell -1} \left( B \left( z_{j_{\ell +1}  ( \dots ( j_{k+1} ( n)) \dots )}^{\ell}, R_{j_{\ell +1}  ( \dots ( j_{k+1} ( n)) \dots )}^{\ell} \right) \right)
\lra m_{\ell  }  \qquad \mbox{ as } n \lra \infty. 
\end{array}
$$
Summing up the above inequalities for $ \ell = 1, \dots , k+1 $,  taking into account (H5) 
and passing to the limit as $ n \lra \infty $ we get 
$$
M \ges \limsup_{ n \ra \infty } \mu_{ j_1( \dots ( j_{k+1} ( n)) \dots )}  (\Xo) \ges  m_1 + \dots + m_{k +1}. 
$$
Hence the properties (H1)$-$(H8) are verified at level $ k+1$ and the induction is complete. 

\medskip

If there is $ k \in \N $ such that $ \al _k = 0$, the process stops after $k$ steps and conclusion (i) in Theorem \ref{main} is satisfied. 
Otherwise we have $ \al _k > 0 $ for all $k$. 
In the latter case the series $ \ds \sum_{ k \ges 0 } \al _k $ converges 
(because  $0 < \al _k < 2 m _{ k +1 } $  and the series $  \ds \sum_{ k \ges 1 } m _k $ converges in view of (H7)).
Therefore $ \al _k \lra 0 $ as $ k \lra \infty$. 

\medskip 

If $ \al _k > 0 $ for all $ k$ we perform a diagonal extraction to get a subsequence  satisfying conclusion (ii) in Theorem \ref{main}. 
Let $ j(n) = j_1 ( \dots ( j_n (n)) \dots ).$
For each $ n \in \N^*$ let 
$$
\begin{array}{l}
B( x_n ^1, r_n ^1 ) = B\left( z_{j_2( \dots ( j_n (n)) \dots)} ^1, R_{j_2( \dots ( j_n (n)) \dots)} ^1  \right)\! ,  \; \; 
\\
\\
B( x_n ^2, r_n ^2 ) = B\left( z_{j_3( \dots ( j_n (n)) \dots) } ^2, R_{j_3( \dots ( j_n (n)) \dots )} ^2 \right) \qquad \mbox{ and so on } \dots, 
\\
\\
B( x_n ^{n-1}, r_n ^{n-1} ) = B \left( z_{j_n (n) }^{ n-1}, R_{j_n (n) }^{ n-1} \right), \; \; 
\\
\\
B( x_n^n, r_n ^n ) = B( z_n^n, R_n ^n ) .
\end{array}
$$
It follows from (H5) that the balls are disjoint. 
For each fixed $k$ the sequence $ ( r_n ^k)_{n \ges 1}$ is increasing and tends to infinity as $ n \lra \infty$, 
and (H6) implies that $ \left( {\mu_{j(n)} } _{\big|  B( x_n ^k,\,  r_n ^k)} \right) _{n \ges k}$ concentrates around $ ( x_n ^k)_{n \ges k }$. 
From (H3) and (H4) we get conclusion (ii) (b). 
By (H2) we have 
$   {\mu_{j(n)} } _{\big|  \Xo \setminus \cup_{i =1}^{\ell} B( x_n ^i, \, r_n ^i)}  = \mu _{j_{\ell +1} ( \dots (j_n (n)) \dots) } ^{\ell}$ if $ \ell < n$ and 
$   {\mu_{j(n)} } _{\big|  \Xo \setminus \cup_{i =1}^n B( x_n ^i, \, r_n ^i)}  = \mu _n ^n$, 
and the concentration functions of these  measures are $ q_{j_{\ell +1} ( \dots (j_n (n)) \dots) } ^{\ell}$ and $ q_n ^n$, respectively. 
Since $ \ds \lim _{t \ra \infty} \left( \limsup_{n \ra \infty }   q_{j_{\ell +1} ( \dots (j_n (n)) \dots) } ^{\ell} (t) \right) \les 
\lim _{t \ra \infty} \left( \limsup_{n \ra \infty }   q_{ n } ^{\ell} (t) \right) = \al _{\ell} $ and $ \al_{\ell } \lra 0 $ as $ \ell \lra \infty$, 
(ii) (c) follows. 
From (H2) and (H8) we get 
$   q_n ^n (n) \les q_{ j_n (n) } ^{ n-1 } (n) < 2 \al _{n-1},  $
hence $  q_n ^n (n)  \lra 0 $ as $ n \lra \infty$.  Then using Lemma \ref{e2} and Lemma \ref{e} (i) we infer that 
$\ds \lim_{ t \ra \infty } \left( \limsup_{ n \ra \infty} {q}_n ^n(t) \right) = 0$ and (ii) (d) is proven.
\hfill
$\Box$

\section{Profile decomposition for bounded sequences in $W^{1, p }(\R^N)$ }
\label{PrSobolev}

The proof of Theorem \ref{Sobolev} relies on Theorem \ref{main} and  the following simple  lemma. 

\begin{Lemma}
\label{L4.1}
Assume that $ 1 \les p < \infty $ and $ q \in (p, p^*)$, where $ p^* = \frac{ Np}{N-p}$ if $ p < N$ and 
$ p^ * = \infty $ if $ p \ges N$. 
There exists $ C> 0$, depending only on $ p, q $ and $N$, such that for any $ u \in W^{1, p } ( \R^N)$ there holds 
$$
\| u \| _{L^q( \R^N)} \les C \left( \sup_{y \in \R^N} \int_{B(y, 1)} |\nabla u |^p + | u |^p \, dx \right)^{ \frac 1p - \frac 1q} \| u \|_{ W^{1, p }(\R^N) } ^{\frac pq}.
$$
\end{Lemma}

{\it Proof. }
 Lemma \ref{L4.1} is known by experts, we give the proof for the sake of completeness. 

By the Sobolev embedding there is $ C_S >0 $ such that for any  ball $ B \subset \R^N$   of radius $1$  and any 
$ w \in W^{1, p }(B)$ there holds
$ \| w \|_{ L^q( B) } \les C_S \| w \|_{W^{1, p }(B)} .$ 
Let $ u \in W^{ 1, p }( \R^N)$. 
For any $ z \in \R^N$ we have 
\beq
\begin{array}{l}
\label{4.1}
\ds \int_{B(z, 1)} | u|^q \, dx \les C_S ^q \left( \int_{B(z, 1)}   |\nabla u |^p + | u |^p \, dx \right) ^{\frac qp} 
\\
\\
\ds \les C_S ^q \left(  \sup_{y \in \R^N} \int_{B(y, 1)}   |\nabla u |^p + | u |^p \, dx \right) ^{\frac qp -1} 
\int_{B(z, 1)}   |\nabla u |^p + | u |^p \, dx .
\end{array}
\eeq
There is an integer $K$ (depending only on $N$) and there is a covering of $ \R^N$ by a family of balls 
$( B( y_n, 1 ))_{n \in \N^*}$ such that each point belongs to at most $K$ balls 
(that is, $ 1 \les \ds \sum_{ n =1}^{\infty} \1_{B(y_n, 1)} \les K$). 
We write (\ref{4.1}) for each ball $B(y_n, 1)$, then sum over $n$ to get the desired conclusion. 
\hfill
$\Box$

\medskip

{\it Proof of Theorem \ref{Sobolev}.}
We consider first the case $ 1 < p < N. $     
By the Sobolev embedding $(u_n)_{n \ges 1}$ is bounded in $ L^{p^*}(\R^N)$. 
Let $ \rho _n = |\nabla u _n |^p + | u_n |^p + |u _n |^{ p^*}, $ so that $ (\rho_n )_{ n \ges 1}$ is bounded in $ L^1( \R^N)$. 
We use Theorem \ref{main} for the sequence of measures $ (\mu_n )_{n \ges 1}$ with densities $ (\rho_n )_{ n \ges 1}$
with respect to the Lebesgue measure. 

If $ (\rho_n )_{n \ges 1}$ is a vanishing sequence it follows from Lemma \ref{L4.1} that $ u_n \lra 0 $ in $ L^q( \R^N)$ for any 
$ q \in (p, p^*)$ anf the conclusion of Theorem \ref{Sobolev} holds with $ V^i = 0 $ for all $i$. 

Assume that $ (\rho_n )_{n \ges 1}$ is not vanishing and consider the mapping $ j : \N^* \lra \N^*$,  the sequences of balls 
$\left (B(x_n ^i, r_n^i) \right)_{n \ges 1}$ and the numbers $ m_i$ given by Theorem \ref{main}. 
If there are only $k$ positive $ m_i$'s  we put $ m_i = 0 $ for $ i > k$.

Let $ v _n ^i = u_{ j(n)} ( \cdot \, + x_n ^i). $

The sequence $ ( v _n ^ 1)_{n \ges 1}$ is bounded in $ W^{1, p }(\R^N)$, 
hence there exists $ \kappa _1 : \N^* \lra \N^*$ increasing and $ V^1 \in W^{1, p }(\R^N) $ such that 
\beq
\label{4.2}
v_{\kappa_1 (n)} ^1 \rightharpoonup V^1 \qquad \mbox{ weakly in } W^{1, p }(\R^N) \mbox{ as } n \lra \infty, 
\eeq
\beq
\label{4.3}
v_{\kappa_1 (n)} ^1 \lra V^1 \quad 
\begin{array}{l} 
\mbox{ strongly in } L^q( B(0, R)) \mbox{ for any } 
 q \in [1, p^*) \mbox{ and any } R > 0 
\\
\mbox{ and a.e. on } \R^N.
\end{array}
 \eeq
Fix $ R > 0$. Since $ \nabla v_{\kappa_1 (n)} ^1 \rightharpoonup \nabla V^1 $ weakly in $ L^p ( B(0, R)) $ and $ r_n ^1 \lra \infty$, we get 
$$
\int_{B(0, R)} |\nabla V^1  |^p \, dx \les \liminf_{ n \ra \infty} \int_{B(0, R)}  |\nabla v_{\kappa_1 (n)}  ^1 |^p \, dx 
\les \liminf_{n \ra \infty} \int_{B(x_{\kappa _1(n)} ^1, r_{\kappa _1(n)} ^1 )}  |\nabla u_{j(\kappa_1 (n)) }   |^p \, dx .
$$
Letting $ R \lra \infty$ we find
\beq
\label{4.4}
\int_{\R^N} |\nabla V^1  |^p \, dx \les 
\liminf_{n \ra \infty} \int_{B(x_{n} ^1, r_{n} ^1 )}  |\nabla u_{j(n) }   |^p \, dx .
\eeq
Using Fatou Lemma and proceeding similarly we discover
\beq
\label{4.5}
\int_{\R^N} | V^1  |^q \, dx \les 
\liminf_{n \ra \infty} \int_{B(x_{n} ^1, r_{n} ^1 )}  | u_{j(n) }   |^q \, dx 
\qquad \mbox{ for all    } q \in [\, p, p^*].
\eeq

Next we show that 
\beq
\label{4.6}
 \1_{B( 0, \, r_{\kappa _1(n)}^1 )}   v_{\kappa_1 (n)} \lra V^1 \qquad \mbox{ strongly in } L^q( \R^N) \mbox{ for any } q \in [\, p, p^*). 
\eeq
Fix $ \e > 0 $. 
Since $ \left( {\mu_{ j(n)}}_{|  B(x_n ^1, r_n^1) } \right) _{n \ges 1}$ concentrates around $  ( x_n)_{ n \ges 1}$, 
there are $ r_{\e } > 0 $ and $ n_{\e } \in \N$ such that 
$ \ds \int _{ B(x _n ^1 , \, r _n ^1) \setminus B( x_n ^1, r_{\e})} \rho_{ j(n)} \, dx < \e$ for all $ n \ges n_{\e}$. 
This implies 
$ \ds \int _{ B(0, \, r _n ^1) \setminus B( 0, r_{\e})}  | v_n ^1 |^p + | v_n ^1 |^{p^*}  \, dx < \e$ for $ n \ges n_{\e } $ 
and using Fatou Lemma we get 
$ \ds \int _{ \R^N  \setminus B( 0, r_{\e})}  | V ^1 |^p + | V ^1 |^{p^*}  \, dx \les \e$. 
By interpolation we find
$\ds \int _{ \R^N  \setminus B( 0, r_{\e})}  | \1_{B(0, r_n ^1)} v_n ^1  |^q \, dx < \e ^{\nu_q} $ and 
$\ds \int _{ \R^N  \setminus B( 0, r_{\e})}  |V ^1 |^q \, dx < \e ^{\nu _q} $
for some $ \nu_q >0$. 
On the other hand, by (\ref{4.3}) there is $ n_{\e}' \ges n_{\e}$ such that 
$\| v_{\kappa_1(n)} ^1 - V ^1 \| _{L^q ( B( 0, r_{\e} ))} < \e $ for $ n \ges n_{\e}'$. 
We conclude that $ \| \1_{B(0, r_n ^1)}  v_{\kappa_1(n) } ^1  - V^1 \|_{L^q ( \R^N)} < \e + 2 \e^{\frac{ \nu_q}{q}}$ 
for all $ n \ges n_{\e}'$
and (\ref{4.6}) is proven.

\medskip

If $ m_2 = 0 $ we take $ V^i = 0 $ for $ i \ges 2. $
If $ m_2 > 0$, proceeding as above we see that there exist an increasing mapping 
$ \kappa_2 : \N^* \lra \N^*$  and $ V^2 \in W^{1, p }(\R^N)$
such that   (\ref{4.2})$-$(\ref{4.6}) hold with 
$ v_{\kappa_1(\kappa_2(n))} ^2 $ and $ V^2$ instead of $ v_{\kappa_1(n)}^1$ and $ V^1$, respectively. 
If $ m_{\ell} > 0 $ for some $ \ell \in \N^*$, by induction we find increasing mappings 
$ \kappa_1, \dots, \kappa_{\ell} : \N^* \lra \N^*$ and 
$ V^1, \dots, V^{\ell } \in W^{1, p }(\R^N)$ such that 
\beq
\label{4.7}
v_{\kappa_1 (\dots ( \kappa_{\ell}( n) ) \dots ) } ^{\ell} \rightharpoonup V^{\ell} \qquad \mbox{ weakly in } W^{1, p }(\R^N) \mbox{ as } n \lra \infty, 
\eeq
\beq
\label{4.8}
v_{\kappa_1 (\dots (\kappa_{\ell} (n)) \dots )} ^{\ell}  \lra V^{\ell} \quad 
\begin{array}{l} 
\mbox{ strongly in } L^q( B(0, R)) \mbox{ for any } 
 q \in [1, p^*) \mbox{ and  } R > 0 
 \\
\mbox{ and  a.e. on } \R^N, 
\end{array}
\eeq
\beq
\label{4.9}
\int_{\R^N} |\nabla V^{\ell}  |^p \, dx \les 
\liminf_{n \ra \infty} \int_{B(x_{n} ^{\ell}, r_{n} ^{\ell} )}  |\nabla u_{j(n) }   |^p \, dx ,
\eeq
\beq
\label{4.10}
\1_{B( 0, \, r_{\kappa _1(\dots ( \kappa_{\ell} (n))\dots ) }^{\ell} )}   v_{\kappa_1 (\dots (\kappa_{\ell}  (n)) \dots) }  ^{\ell} \lra V^{\ell} \qquad \mbox{ strongly in } L^q( \R^N) \mbox{ for any } q \in [\, p, p^*). 
\eeq

If there is $ \ell_0 $ such that  $ m_{\ell_0 } > 0$ and $ m_i = 0 $ for $ i > \ell_0$, we put $ V^i = 0 $ for $ i > \ell_0$ and 
$ \kappa = \kappa_1 \circ \dots \circ \kappa_{\ell_0}$. 
Otherwise we find $ V^i $ as above for all $ i \in \N^*$ and we put 
$ \kappa ( n )= \kappa_1( \dots ( \kappa_n (n)) \dots ).$ 
We show that the conclusion of Theorem \ref{Sobolev} holds true for the subsequence $\left( u_{j(\kappa(n))}\right)_{n \ges 1}$ 
and the functions $V^i $ as above.

Fix $ \chi \in C_c ^{\infty}(\R^N)$ such that $ 0 \les \chi \les 1$, $ \chi = 1$ on $B(0, \frac 12)$ and $ \mbox{supp}( \chi)\subset B(0,1)$. 
For $ n \ges k $ the balls $ B( x_n ^1, r_n ^1), \dots B( x_n ^k, r_n ^k)$ are disjoint and it is easy to see that 
\beq
\label{4.11}
1= \chi \left( \frac{ \cdot - x_n ^1}{ r_n ^1} \right) + \dots + \chi \left( \frac{ \cdot - x_n ^k}{ r_n ^k} \right) 
+ \prod_{i = 1}^k  \left( 1 - \chi \left( \frac{ \cdot - x_n ^i}{ r_n ^i} \right) \right).
\eeq

Fix $ k \in \N^*$ such that $ m_k > 0 $. For $ n \ges k $ we may write 
\beq
\label{4.12}
\begin{array}{rcl}
u_{j ( \kappa (n))} & = & \ds  \sum_{i =1}^k V ^i ( \cdot - x_{\kappa(n)} ^i ) + 
\sum_{i =1}^k \left( \chi \left( \frac{ \cdot - x_{\kappa(n)} ^i  }{r_{\kappa(n)} ^i  } \right) u_{j( \kappa(n))} -  V ^i ( \cdot - x_{\kappa(n)} ^i ) \right)
\\
\\
& & \ds
+ \left( \prod _{i =1}^k \left( 1 - \chi \left( \frac{ \cdot - x_{\kappa(n)} ^i  }{r_{\kappa(n)} ^i  } \right) \right) \right) u_{j( \kappa(n))}.
\\
\\
& = &  \ds  \sum_{i =1}^k V ^i ( \cdot - x_{\kappa(n)} ^i ) + v_{n, 1}^k + v_{n, 2}^k. 
\end{array}
\eeq
Let $ w_n ^k = v_{n, 1}^k + v_{n, 2}^k. $
The fact that the sequence $ \left( \1_{B( x_{ \kappa (n)} ^i, \,  r_{ \kappa (n)}^i) }  \rho_{j ( \kappa (n))} \right)_{n \ges 1}$ 
concentrates around $ \left(x_{\kappa (n) }^i \right)_{n \ges 1}$ and (\ref{4.10}) imply that 
$ \| v_{n, 1}^k \| _{L^q  (\R^N) } \lra 0 $ as $ n \lra \infty$ for any $ q \in [p , p^*)$. 
It is  easy to see that 
$$
|  v_{n, 2}^k | ^p + |\nabla  v_{n, 2}^k |^p \les \rho_{j(\kappa(n))} \1_{ \R^N \setminus \cup_{i=1}^k B( x_{j ( \kappa (n))} ^i,\,  r_{j ( \kappa (n))}^i)} + f_n ^k , 
$$ 
where $ \| f_n ^k \|_{L^1( \R^N)} \lra 0 $ as $n \ra \infty$. 
Theorem \ref{main} implies that 
$$
\ds \limsup_{n \ra \infty} \left( \sup_ { z \in \R^N} \int_{B(z, 1)} |  v_{n, 2}^k |^p + |\nabla  v_{n, 2}^k |^p \, dx \right) \lra 0 \qquad \mbox{  as  } k \lra \infty.
$$ 
Then using Lemma \ref{L4.1} we get $ \ds \lim_{k \ra \infty} \left( \limsup_{n \ra \infty} \|  v_{n, 2}^k \|_{L^q ( \R^N )} \right) = 0$ 
for any $ q \in ( p,  p^*)$
and conclusion  (ii) in Theorem  \ref{Sobolev} follows. 

\medskip

Using (\ref{4.10}), the definition of $ v_{n, 2}^k $ (see (\ref{4.12})) and the fact that 
$ \| v_{n, 1}^k \| _{L^p  (\R^N) } \lra 0 $ as $ n \lra \infty$  we get 
$$
\begin{array}{l}
\ds \int_{\R^N} |u _{j ( \kappa (n))} |^p \, dx 
= \ds \sum_{i = 1}^k \int_{B( x_{ \kappa (n)} ^i,\,  r_{\kappa (n) }^i) } |u _{j ( \kappa (n))} |^p \, dx  + 
\int_{\R^N \setminus \cup_{i=1}^k B( x_{  \kappa (n)} ^i,\,  r_{ \kappa (n) }^i)} |u _{j ( \kappa (n))} |^p \, dx
\\
\\
= \ds \sum_{i = 1}^k  \|V^i \| _{L^p(\R^N) }^p + \|  v_{n, 2}^k \| _{L^p(\R^N) }^p + o(1) 
\\
\\
= \ds \sum_{i = 1}^k  \|V^i \| _{L^p(\R^N) }^p + \|  w_{n}^k \| _{L^p(\R^N) }^p + o(1) 
\qquad \mbox{ as } n \lra \infty 
\end{array}
$$
and (iii) is proven.

\medskip

Let $ h_n ^i (x) = \ds \chi \left(\frac{x}{r_{\kappa(n)} ^i } \right) u_{j( \kappa (n))} ( \cdot + x_{\kappa(n)} ^i ).$
For all $n$ such that $r_{\kappa(n)} ^i \ges 1$ we have  
\beq
\label{4.13}
\begin{array}{l}
\ds \int_{\R^N} \Big|  |\nabla  u_{j( \kappa (n))} |^p - |\nabla v_{n, 2} ^k |^p - \sum_{i =1}^k | \nabla h_n ^i ( x - x_{\kappa(n)} ^i ) |^p \Big| dx 
\\
\\
\ds \les C \sum_{i=1}^k \int_{B(  x_{\kappa(n)} ^i , r_{\kappa(n)} ^i) \setminus B(  x_{\kappa(n)} ^i , \frac 12 r_{\kappa(n)} ^i) } |\nabla  u_{j( \kappa (n))} |^p + | u_{j( \kappa (n))} |^p \, dx. 
\end{array}
\eeq
and the right hand side in the above inequality tends to zero as $ n \lra \infty $ because the sequence 
$ \left( \1_{B(  x_{\kappa(n)} ^i , \, r_{\kappa(n)} ^i)}   \rho_{j(\kappa (n))} \right)_{n \ges 1} $ 
concentrates around $( x_{\kappa(n)} ^i )_{n \ges 1}$. 

Using (\ref{4.7}) it is easy to see that $ h_n ^i \rightharpoonup V^i$ weakly in $ W^{1, p }( \R^N)$. 
Assume that $ p = 2$. 
Then we have 
\beq
\label{4.14}
\begin{array}{l}
\ds \int_{\R^N} |\nabla h_n ^i ( x - x_{\kappa(n)} ^i ) |^2 \, dx = \int_{\R^N} |\nabla h_n ^i ( x ) |^2 \, dx  
\\
\\
= \| \nabla V ^i \|_{L^2( \R^N)}^2 + \| \nabla h_n ^i - \nabla V ^i \|_{L^2( \R^N)}^2 + 2 \langle \nabla h_n ^i - \nabla V ^i , \nabla V^i \rangle _{L^2}
\\
\\
= \| \nabla V ^i \|_{L^2( \R^N)}^2 + \| \nabla h_n ^i - \nabla V ^i \|_{L^2( \R^N)}^2 + o(1). 
\end{array}
\eeq
Since $ r_{\kappa(n)} ^i \lra \infty$ for all $i$ and $| x_{\kappa(n)} ^i - x_{\kappa(n)} ^j | \lra \infty$ as $ n \lra \infty$, 
it is straightforward to check that 
\beq
\label{4.15}
 \| \nabla v_{n,1} ^k \|_{L^2( \R^N)}^2  = \sum_{i =1}^k  \| \nabla h_n ^i - \nabla V ^i \|_{L^2( \R^N)}^2 +  o(1) \qquad \mbox{ and } 
\eeq
\beq
\label{4.16}
\| \nabla w_n ^k \|_{L^2( \R^N)}^2 = \| \nabla v_{n,1} ^k \|_{L^2( \R^N)}^2 + \| \nabla v_{n,2} ^k \|_{L^2( \R^N)}^2 + o(1) 
\eeq
as $ n \lra \infty$. 
Then (iv) follows from (\ref{4.13})$-$(\ref{4.16}).

\medskip

Assume next that $ N \les p < \infty.$
Fix an increasing sequence $ (p_{\ell})_{\ell \ges 1} \subset (p, \infty)$, $p_{\ell} \lra \infty. $
By the Sobolev embedding for any $ \ell $ there is $ C_{\ell} > 0 $ such that for any $ u \in W^{1, p }(\R^N)$ there holds
$ \| u \|_{L^{p _{\ell}} ( \R^N )} \les C_{\ell}  \| u\|_{W^{1, p }(\R^N)} $. 
Let $ A = \ds \sup_{n \ges 1} \| u_n \|_{W^{1, p }(\R^N)}. $
Let 
$$
\rho_n = |\nabla u_n |^p + | u_n |^p + \sum_{\ell = 1}^{\infty} \frac{1}{2^{\ell} C_{\ell}^{p_{\ell}} A^{p_{\ell}} } | u_n |^{p_{\ell}} .
$$
Then $(\rho_n)_{n \ges 1}$ is bounded in $ L^1( \R^N)$ and for any $ s \in [\, p , \infty)$ there exist $ a_s>0$, $ \nu_s >0$ such that 
for all $n$ and all measurable sets $E$, 
$$
\int_{E} | u_n |^s \, dx \les a_s \left(  \int_{E}  \rho_n \, dx \right)^{ \nu _s}. 
$$
We apply Theorem \ref{main} to 
$(\rho_n)_{n \ges 1}$. 
The rest of the proof is as in the case $ p \in (1, N).$
\hfill
$\Box$

\bigskip

\noindent
{\bf Acknowledgements. } 
I am indebted to Patrick G\'erard for suggesting many important references.
I am very grateful to Franck Barthe, David Chiron,  
 Radu Ignat, Stefan Le Coz, Michel Ledoux, Yvan Martel, Frank Merle  and Nicola Visciglia for interesting and helpful discussions.

\end{document}